\documentclass[12pt]{amsart}
\usepackage{hyperref}

\usepackage[utf8]{inputenc}
\usepackage[T1]{fontenc}

\usepackage[]{amsfonts}
\usepackage{amssymb}
\usepackage{amsmath}
\def\couleur(#1 #2 #3)
	{
\special{ps::[end]}}

\def\bx#1{\setbox1=\hbox{\kern3pt{#1}\kern3pt}			
 \dimen1=\ht1 \advance\dimen1 by 3pt \dimen2=\dp1 \advance\dimen2 by 3pt
 \setbox1=\hbox{\vrule height\dimen1 depth\dimen2\box1\vrule}%
 \setbox1=\vbox{\hrule\box1\hrule}%
 \advance\dimen1 by .4pt \ht1=\dimen1
 \advance\dimen2 by .4pt \dp1=\dimen2 \box1\relax}

\def\wbb#1{\kern#1em}
\def\vci{\vrule  width.02em height1.47ex depth-.0ex}		
\def\11{{\rm\wbb{.2}\vci\wbb{-.37}1}}

\def\underset#1#2{\mathrel{\mathop{\kern0pt #2}\limits_{#1}}}

\def\overset#1#2{\mathrel{\mathop{\kern0pt #2}\limits^{#1}}}

\newtheorem{Thrm}{Theorem}[section]
\newtheorem{Lmm}[Thrm]{Lemma}
\newtheorem{Dfnt}[Thrm]{Definition}
\newtheorem{Prps}[Thrm]{Proposition}
\newtheorem{Rmrq}[Thrm]{Remark}
\newtheorem{Exmp}[Thrm]{Example}

\begin{document}

\title{Subspaces of $\displaystyle H^{p}$ linearly homeomorphic to $l^{p}.$}

\author{Amar E., Chalendar I., Chevreau B.}
\maketitle
 \ \par 
\ \par 
\renewcommand{\abstractname}{Abstract}

\begin{abstract}
We present two fast constructions of weak*-copies of $\ell ^{\infty
 }$  in $H^{\infty }$ and show that such copies are necessarily
 weak*-complemented. Moreover, via a Paley-Wiener type of stability
 theorem for bases, a connection can be made in some cases between
 the two types of construction, via interpolating sequences (in
 fact these are at the basis of the second construction). Our
 approach has natural generalizations where $H^{\infty }$ is
 replaced by an arbitrary dual space and $\ell ^{\infty }$ by
 $\ell ^{p}$ ($1\leq p\leq \infty $) relying on the notions of
 generalized interpolating sequence and bounded linear extension.
 An old (very simple but unpublished so far) construction of
 bases which are Besselian but not Hilbertian finds a natural
 place in this development.\ \par 
\end{abstract}
\ \par 

\section{Introduction.}
\quad This note is inspired by the construction in~\cite{BiF6}  of
 special weak*-closed subspaces in the algebra $H^{\infty }$
 (of bounded analytic functions in the open unit disc ${\mathbb{D}}$).
 The last two named authors realized that it could be used to
 obtain weak* homeomorphic copies of $\ell ^{\infty }$ in $H^{\infty
 }$ in a very elementary fashion (at the expense of "forgetting"
 several additional technicalities needed by the author for later
 purposes). Then, the first author observed that another fast
 construction could be achieved via interpolating sequences using
 a deep result of P. Beurling. Thus below (Section~\ref{F5})
 we present both types of construction. In fact a (relatively
 easy) characterization of weak*-continuous linear mappings from
 $\ell ^{\infty }$ to $H^{\infty }$(Prop.~\ref{FB0}) enables
 us to give a convenient description of "weak*-copies" of $\ell
 ^{\infty }$  in $H^{\infty }$ and, in particular, showing further
 on (Theorem~\ref{FB2})  that any such copy is weak*-complemented
 in $H^{\infty }.$\ \par 
\quad The starting point of our first construction is a countable separated
 set (i.e. consisting of isolated points) $S:=\lbrace a_{n};\
 n\in {\mathbb{N}}\rbrace $  in the unit circle ${\mathbb{T}}.$
 When the closure of this set has arc-length measure $0$ the
 Rudin-Carleson interpolating theorem (see Remark 2 following
 Theorem~\ref{F6}) enables us to prove the existence of a weak*-closed
 subspace ${\mathcal{E}}$ of $H^{\infty }$ whose predual has
 an $\ell ^{1}-$ basis (see below Section~\ref{FB1} for the terminology)
 consisting of evaluations at the points $a_{n}.$ Then, a straightforward
 application of a Paley-Wiener type of stability theorem for
 bases (in the predual of ${\mathcal{E}}$) leads to the fact
 that any sequence $(b_{n})_{n}$ in ${\mathbb{D}}$ with $b_{n}$
  "close enough" to $a_{n},$ is an interpolating sequence (i.e.
 the map $h\rightarrow (h(b_{n}))_{n}$ from $H^{\infty }$ is
 onto $\ell ^{\infty }$) (see Theorem~\ref{F4} and Theorem~\ref{F7}).
 In fact, this profusion of interpolating sequences in the "vicinity"
 of $(a_{n})_{n}$ can be obtained directly from Carleson characterization
 of such sequences under the mere separation of the sequence
 $(a_{n})_{n}\ ;$ moreover an idea from~\cite{BiF8} leads to
 an explicit construction of a weak*-copy of $\ell ^{\infty }$
  for any of these interpolating sequences (cf. Theorem~\ref{F8}).\ \par 
\quad The above theme leads quite naturally to the question (already
 considered in the literature) of how to get copies of $\ell
 ^{p}$ in the Hardy space $H^{p}$ for $1\leq p<\infty $ ; this
 question is related to the notion of Besselian and Hilbertian
 sequences and we briefly discuss it in Section~\ref{F9}. Along
 the way we give an elementary construction (coming from the
 first named author's thesis but never published elsewhere) of
 a Hilbertian basis which is not Besselian (the first such known
 example being due to Babenko~\cite{BiF9}).\ \par 
\quad During the elaboration of this paper we discovered that certain
 of our results were already known, e.g., the construction of
 copies of $\ell ^{\infty }$ via interpolating sequences is in
 Chap. VII of~\cite{BiF10} or particular cases of more general
 ones (e.g., Theorem~\ref{FB2} which is a corollary of a deep
 result of~\cite{BiF11}). However, in line with the above mentioned
 elementary character of the first construction (and starting
 point of this work), we have chosen to make this note as self-contained
 as possible, thus illustrating in (we hope) an easy fashion
 the interplay between some pieces of Banach space geometry and
 function theory.\ \par 
\quad Some terminology and preliminary material are developed in Section~\ref{FB1}.\
 \par 

\section{Notations and Preliminaries.~\label{FB1}}
\quad As usual, $A({\mathbb{D}})$ (${\mathbb{D}}$ open unit disk) denotes
 the disk algebra consisting of those functions continuous on
 the closed unit disk and analytic in ${\mathbb{D}}.$ For $1\leq
 p\leq \infty $  the Hardy space $H^{p}$ can be defined as the
 norm-closure (weak*-closure if $p=\infty $) of $A({\mathbb{D}})$
 in $L^{p}(=L^{p}(m))$ (where $m$ denotes the normalized Lebesgue
 measure on the unit circle ${\mathbb{T}}$). Via Fatou's lemma
 the space $H^{p}$ can (and will) be identified as the space
 of analytic functions $f$ in ${\mathbb{D}}$ satisfying\ \par 
\quad \quad \quad $\displaystyle ({\left\Vert{f}\right\Vert}_{p}:=)\mathrm{s}\mathrm{u}\mathrm{p}_{0\leq
 r<1}{\left\Vert{f_{r}}\right\Vert}_{p}<\infty $\ \par 
where $f_{r}$ is defined on ${\mathbb{T}}$ by $f_{r}(w)=f(rw)$
 and the above weak*-closure refers to the duality $L^{\infty
 }(=L^{\infty }(m))=(L^{1})^{*}.$\ \par 
\quad In particular, $H^{\infty }$ can be seen as the Banach algebra
 of bounded analytic functions in ${\mathbb{D}}$ (equipped with
 the "sup" norm: $\ {\left\Vert{h}\right\Vert}_{\infty }:=\mathrm{s}\mathrm{u}\mathrm{p}_{z\in
 {\mathbb{D}}}\left\vert{h(z)}\right\vert $)  and, since $H_{0}^{1}$
 is the preannihilator of $H^{\infty }$ in the above duality,
 $H^{\infty }$ is also the dual of the quotient space $L^{1}/H_{0}^{1}.$
 If $\lbrack f\rbrack $ denotes the equivalence class of $f\in
 L^{1}$ and $h\in H^{\infty }$ we have thus\ \par 
\quad \quad \quad $\displaystyle \ {\left\langle{\lbrack f\rbrack ,h}\right\rangle}={\left\langle{f,h}\right\rangle}=\int{fh\
 (:=}\int_{{\mathbb{T}}}^{\ }{fhdm=}\frac{1}{2\pi }\int_{0}^{2\pi
 }{f(e^{it})h(e^{it})dt)}.$\ \par 
\quad Among the elements in the predual $L^{1}/H_{0}^{1}$ of $H^{\infty
 }$  we have the evaluations $E_{\lambda }$ at points $\lambda
 \in {\mathbb{D}}$  ($\ {\left\langle{E_{\lambda },h}\right\rangle}=h(\lambda
 )$) ; in fact $E_{\lambda }=\lbrack P_{\lambda }\rbrack $  where
 $P_{\lambda }$ is the Poisson kernel at $\lambda $ ($P_{\lambda
 }(z)=\frac{1-\left\vert{\lambda }\right\vert ^{2}}{\left\vert{1-\bar
 \lambda z}\right\vert ^{2}}$). Of course these evaluations define
 (by restriction) weak*-continuous linear functionals on any
 weak*-closed subspace ${\mathcal{E}}$ of $H^{\infty }$ and will
 thus be identified (keeping the same notation $E_{\lambda }$)
 with elements in the predual ${\mathcal{Q}}_{{\mathcal{E}}}$
  of ${\mathcal{E}}$ (this predual ${\mathcal{Q}}_{{\mathcal{E}}}$
 can be seen either as a quotient of $L^{1}$ or of $L^{1}/H_{0}^{1}$).
 In some instances (see again Remark 2 after Theorem~\ref{F6})
 evaluations at points of ${\mathbb{T}}$ make sense and are weak*-continuous.\
 \par 
\quad We recall that a sequence $(h_{n})$ in $H^{\infty }$ converges
 weak* to $0$ if and only if it is bounded and converges pointwise
 to $0$ on ${\mathbb{D}}.$\ \par 
\quad We will also need some general duality facts. Given two Banach
 spaces $X,Y,\ {\mathcal{L}}(X,Y)$ will denote the Banach space
 of bounded linear operators from $X$ to $Y.$ We recall the following
 fundamental and classical result (cf.~\cite{BiF5}, Theorem 4.14.).\ \par 

\begin{Thrm}
For $S\in {\mathcal{L}}(X,Y)$ the following three assertions
 are equivalent : \par 
a) $S$ has closed range;\par 
b) $S^{*}$ has (norm-) closed range;\par 
c) $S^{*}$ has weak*-closed range.
\end{Thrm}
\quad As an almost immediate corollary we have the following useful
 proposition.\ \par 

\begin{Prps}
~\label{F10}Let $X,Y$ be Banach spaces and ${\mathcal{M}}$ a
 (linear) subspace of $X^{*}.$\par 
a) Suppose ${\mathcal{M}}$ is weak*-homeomorphic to $Y^{*}$ via
 a (linear) map $T$ then ${\mathcal{M}}$ is weak*-closed (and
 $T$  is norm-invertible).\par 
b) Conversely if a linear map $T\ :Y^{*}\rightarrow X^{*}$ is
 weak*-continuous and bounded below then it implements a weak*-homeorphism
 between $Y^{*}$ and the (weak* as well as norm-closed) subspace
 ${\mathcal{M}}=T(Y^{*}).$
\end{Prps}
\quad Proof\ \par 
Everything is an immediate consequence of the classical duality
 result recalled above except the first part of a) for which
 (by virtue of this same result) we just need to show that $T$
 is bounded below. Suppose not: then there exists an unbounded
 sequence $(y_{n})_{n}$  in $Y^{*}$ such that $\ {\left\Vert{Ty_{n}}\right\Vert}\rightarrow
 0$ ; consequently the sequence $Ty_{n}$ converges weak* to $0$
 and so does the sequence $(y_{n})_{n}$ in contradiction with
 its unboundedness. $\displaystyle \hfill\blacksquare $\ \par 
\ \par 
\quad Spaces of complex-valued sequences will intervene throughout
 this paper. We denote by $\ell ^{0}$ (resp. $c_{00}$) the linear
 space of all complex sequences (resp. finitely supported sequences)
 and by $(\epsilon _{n})_{n\in {\mathbb{N}}},$ ($\epsilon _{n}=(\delta
 _{n,j})_{j\in {\mathbb{N}}}$)  the canonical (algebraic) basis
 of $c_{00}.$ The Banach spaces $c_{o},\ \ell ^{p},\ 1\leq p\leq
 \infty $  have their usual meaning.\ \par 
\quad We recall that a Schauder basis in a Banach space $X$ (resp.
 a w*-Schauder basis in a dual space $X=({\mathcal{Q}}_{X})^{*}$)
 is a sequence $(x_{n})_{n}$  in $X$ such that any $x\in X$ has
 a unique norm expansion $x=\Sigma _{n}\alpha _{n}x_{n}$  (resp.
 weak*-convergent expansion $x=\Sigma _{n}^{w*}\alpha _{n}x_{n},$
 that is, for any $L\in {\mathcal{Q}}_{X}$ we have $\ {\left\langle{L,x}\right\rangle}=\sum_{n}{\alpha
 _{n}{\left\langle{L,x_{n}}\right\rangle}}$)  with $(\alpha _{n})_{n}\in
 \ell ^{0}.$ Of course $(\epsilon _{n})_{n}$  is a Schauder basis
 for $c_{0}$ and $\ell ^{p}$ (for $p<\infty $)  and a weak*-Schauder
 basis for $\ell ^{\infty }.$ (The reader is referred to~\cite{BiF12},~\cite{BiF13}
  for more details on bases.)\ \par 
\quad We recall also that two Schauder bases (resp. two w*-Schauder
 bases)  $(x_{n})_{n}$ in $X$ and $(y_{n})_{n}$ in $Y$ are said
 to be equivalent (resp. w*-equivalent) if there exists $A\in
 {\mathcal{L}}(X,Y)$ invertible (resp. w*-homeomorphism) such
 that $Ax_{n}=y_{n}$ for all $n.$ We will say that a Schauder
 basis in $X$ is a $c_{0}-$ basis (resp. $\ell ^{p}-$ basis)
 if it is equivalent to the basis $(\epsilon _{n})_{n}$  (in
 the corresponding space and, of course, equivalence meaning
 w*-equivalence in the case $p=\infty $). Clearly $X$ has an
 $\ell ^{1}-$ basis if and only if $X^{*}$ has an $\ell ^{\infty
 }-$ basis. Thus our initial goal can be reformulated as finding
 weak*-closed subspaces with an $\ell ^{\infty }-$ basis (or
 whose predual has an $\ell ^{1}-$ basis).\ \par 
\quad A first step in that direction is the following useful characterization
 of weak*-continuous linear maps from $\ell ^{\infty }$ into
 $H^{\infty }.$\ \par 

\begin{Prps}
~\label{FB0}  Let $(g_{n})_{n\in {\mathbb{N}}}$ be a sequence
 of functions in $H^{\infty }.$ The following 3 conditions are equivalent:\par 
\quad (i) There exists a weak*-continuous linear map $T\ :\ \ell ^{\infty
 }\rightarrow H^{\infty }$ such that $T\epsilon _{j}=g_{j}\ ,\
 j\in {\mathbb{N}}.$\par 
\quad (ii)             $M:=\mathrm{s}\mathrm{u}\mathrm{p}_{z\in {\mathbb{D}}}\sum_{n\in
 {\mathbb{N}}}{\left\vert{g_{n}(z)}\right\vert }<\infty .$\par 
\quad (iii) There exists a Borel subset $B$ of ${\mathbb{T}}$ of full
 measure (i.e. $m(B)=1$) such that\par 
$\ \ \ \ \ \ \ \ \ \ \ \ \tilde M:=\mathrm{s}\mathrm{u}\mathrm{p}_{z\in
 {\mathbb{D}}}\sum_{n\in {\mathbb{N}}}{\left\vert{g_{n}(z)}\right\vert
 }<\infty .$\par 
\quad When these conditions are fullfilled we have $M=\tilde M$  and
 $T=S^{*}$ with $S$ (from $L^{1}/H_{0}^{1}$ into $\ell ^{1}$)
 defined by $\lbrack f\rbrack \rightarrow (\int_{{\mathbb{T}}}{fg_{n})_{n\in
 {\mathbb{N}}}}$ and $M=\tilde M={\left\Vert{S}\right\Vert}(={\left\Vert{T}\right\Vert}).$\par
 
Furthermore, in case the $g_{n}$'s all belong to the disk algebra,
 they are equivalent to\par 
\quad (ii)'      $M_{1}:=\mathrm{s}\mathrm{u}\mathrm{p}_{z\in {\mathbb{D}}}\sum_{n\in
 {\mathbb{N}}}{\left\vert{g_{n}(z)}\right\vert }<\infty $ (and $M_{1}=M$ ).
\end{Prps}
\quad Proof\ \par 
The last assertion (equivalence of, say, (ii) and (ii)' under
 membership of the $g_{n}$'s in the disk algebra) is obvious
 and stated here just for "the record".\ \par 
\quad (i) $\Rightarrow $ (ii)\ \par 
Let $T$ be a weak*-continuous linear map from $\ell ^{\infty
 }$ in $H^{\infty }$ such that $T\epsilon _{j}=g_{j},\ j\in {\mathbb{N}}.$
 It is straightforward to check that (as indicated at the (iii)
 of the proposition~\ref{FB0}) $T$ is the adjoint of the (bounded
 linear) map $S$ (from $L^{1}/H_{0}^{1}$ into $\ell ^{1}$) defined
 by $\lbrack f\rbrack \rightarrow (\int_{{\mathbb{T}}}{fg_{n})_{n\in
 {\mathbb{N}}}}$ ; thus, for $\lambda \in {\mathbb{D}},\ S(E_{\lambda
 })=(\int{P}_{\lambda }g_{n})_{n}=(g_{n}(\lambda ))_{n}$  and
 we have $\ \sum_{n\in {\mathbb{N}}}{\left\vert{g_{n}(\lambda
 )}\right\vert }(={\left\Vert{S(E_{\lambda })}\right\Vert}_{1})\leq
 {\left\Vert{S}\right\Vert}$ (since the -characters- $E_{\lambda
 }$ are of norm $1$). Hence we have (ii) with $M\leq {\left\Vert{S}\right\Vert}.$\
 \par 
\quad (ii) $\Rightarrow $ (iii)\ \par 
For each $n$ there exists a subset $\omega _{n}$ of ${\mathbb{T}},$
 of measure $0,$ such that at any $w\in {\mathbb{T}}\backslash
 \omega _{n},\ g_{n}(w)$ is the limit of $g_{n}(z)$ as $z$ approaches
 $w$ nontangentially, $g_{n}(w)=NT-\mathrm{l}\mathrm{i}\mathrm{m}_{z\rightarrow
 w,\ z\in {\mathbb{D}}}g_{n}(w).$ The set $\omega :=\bigcup _{n}\omega
 _{n}$ is also of measure $0$ and for $z\in B:={\mathbb{T}}\backslash
 \omega ,\ n\in {\mathbb{N}}$  we have\ \par 
\quad \quad \quad $\displaystyle \ \sum_{k\leq n}{\left\vert{g_{k}(w)}\right\vert
 }=NT-\lim \ _{z\rightarrow w,\ z\in {\mathbb{D}}}\sum_{k\leq
 n}{\left\vert{g_{k}(z)}\right\vert }\leq M.$\ \par 
\quad Therefore the sequence $(g_{n})_{n\in {\mathbb{N}}}$  satisfies
 (iii) with $\tilde M\leq M.$\ \par 
\quad (iii) $\Rightarrow $ (i)\ \par 
Let $c=(c_{n})_{n}\in \ell ^{\infty }$ ; the series $\ \sum_{n}{c_{n}g_{n}}$
 is absolutely convergent on $B$ (indeed, for $z\in B,\ \sum_{n}{\left\vert{c_{n}g_{n}(z)}\right\vert
 }\leq \tilde M{\left\Vert{c}\right\Vert}_{\infty }$). Hence,
 the sequence of partial sums ($G_{N}=\sum_{n\leq N}{c_{n}g_{n}}$)
  is uniformly bounded and, via an application of the Lebesgue
 dominated convergence theorem to the equalities $\ {\left\langle{f,G_{N}}\right\rangle}=\int{fG_{N}},\
 f\in L^{1},\ N\in {\mathbb{N}}$ we obtain its weak*-convergence
 in $L^{\infty }.$ Since $H^{\infty }$ is weak*-closed in $L^{\infty
 }$ the limit $Tc:=\sum_{n}{c_{n}g_{n}}$ belongs to $H^{\infty
 }.$ Moreover the (linear) map $S_{o}:\lbrack f\rbrack \rightarrow
 ({\left\langle{f,g_{n}}\right\rangle})_{n}=(\int{fg_{n}})_{n}$
  from $L^{1}/H_{0}^{1}$ in (a priori) $\ell ^{0}$ is in fact
 $\ell ^{1}$ valued and bounded ; indeed\ \par 
\quad \quad \quad $\displaystyle \ \sum_{n}{\left\vert{\int{fg_{n}}}\right\vert
 }\leq \sum_{n}{\int{\left\vert{fg_{n}}\right\vert }}\leq \int{\left\vert{f}\right\vert
 \sum_{n}{\left\vert{g_{n}}\right\vert }}\leq \tilde M{\left\Vert{f}\right\Vert}_{1}.$\
 \par 
\quad The duality relation\ \par 
\quad \quad \quad $\displaystyle \ {\left\langle{Sf,c}\right\rangle}={\left\langle{f,Tc}\right\rangle}\
 (f\in L^{1},\ c\in \ell ^{\infty })$\ \par 
yields the equality $T=S^{*}$ and the weak*-continuity of $T.$
 Note also that $\ {\left\Vert{T}\right\Vert}={\left\Vert{S}\right\Vert}\leq
 \tilde M;$  combined with the above inequalities $M\leq {\left\Vert{S}\right\Vert}$
 and $\tilde M\leq M$ this proves the equality $M=\tilde M$ and
 concludes the proof. $\displaystyle \hfill\blacksquare $\ \par 
\ \par 
\quad Thus, having a sequence $(g_{n})_{n}$ satisfying Assertion (ii)
 in the above proposition we need "something more" to ensure
 the boundedness below of $T.$ The idea is to have $\ \left\vert{g_{n}}\right\vert
 $ equal or close to $1$ at some point $a_{n}$ or on some subset
 (these points or subsets being sufficiently "separated" for
 different values of $n$) and small elsewhere, $0$ if possible
 at the $a_{k}$ 's for $k\neq n$ (note that this idea clearly
 points the way towards interpolating sequences).\ \par 

\section{Constructions of weak*-copies of $l^{\infty }$ in $\displaystyle
 H^{\infty }.$~\label{F5}}

\subsection{An elementary construction.~\label{F12}}
\ \par 
\quad It is based on the following simple lemma. Here $\bar {\mathbb{D}}$
  is the closed unit disk and for a subarc $\Gamma $ of ${\mathbb{T}}$
  with (distinct) endpoints $w_{1},\ w_{2}$ we denote by $\Delta
 _{\Gamma }$  the union of $\Gamma $ with the open domain (contained
 in ${\mathbb{D}}$)  whose boundary consists of $\Gamma $ together
 with the closed segment of endpoints $w_{1},\ w_{2}$ (in case
 $\Gamma $ is a closed arc we exclude $w_{1},\ w_{2}$ from $\Delta
 _{\Gamma }$)\!\!\!\! .\ \par 

\begin{Lmm}
Let $\Gamma $ be a proper subarc of  ${\mathbb{T}}$ whose complement
 is not a singleton, $a$ a point of $\Gamma $ and $\delta >0,$
 then there exists a function $g\in A({\mathbb{D}})$ such that\par 
\quad (i) $g(a)=1$ and $\ \left\vert{g(z)}\right\vert <1$ for $z\in
 \bar {\mathbb{D}}\backslash \lbrace a\rbrace ,$ and\par 
\quad (ii) $\ \left\vert{g}\right\vert <\delta $ on $\bar {\mathbb{D}}\backslash
 \Delta _{\Gamma }.$
\end{Lmm}
\quad Proof\ \par 
Take $g(z)=(\frac{\bar az+1}{2})^{q}$ (or $(\frac{1}{2-\bar az})^{q}$
 )  with $q$ sufficiently large. $\displaystyle \hfill\blacksquare $\ \par 
\ \par 
\quad We now give ourselves a sequence of open subarcs of ${\mathbb{T}},\
 (\Gamma _{n})_{n\in {\mathbb{N}}}$ whose closures (in ${\mathbb{T}}$)
  are pairwise disjoint, a sequence $(a_{n})_{n}$ such that,
 for each $n,\ a_{n}\in \Gamma _{n}$ and strictly positive scalars
 $\epsilon ,\ \delta _{n}$  such that $\ \sum_{n\in {\mathbb{N}}}{\delta
 _{n}}<\epsilon <1.$ (We note that if we wish to start this procedure
 with a given sequence $(a_{n})_{n}$ the only condition required
 on this sequence is that it be separated.)\ \par 
\quad We claim that the sequence given by application of the lemma
 (that is, $(g_{n})_{n}$ in the disk algebra $A({\mathbb{D}})$
 such that, for each $n,\ g_{n}(a_{n})=1,\ \left\vert{g_{n}}\right\vert
 <1$ on ${\mathbb{T}}\backslash \lbrace a_{n}\rbrace ,$ and $\
 \left\vert{g_{n}}\right\vert <\delta _{n}$ on ${\mathbb{T}}\backslash
 \Gamma _{n}$)  satisfies Assertion (i) of the proposition~\ref{FB0}
 with $M=1+\epsilon .$\ \par 
\quad Indeed for $z\in \Gamma _{n}$ we have $\ \left\vert{g_{n}(z)}\right\vert
 \leq 1$ and for $k\neq n$ $\ \left\vert{g_{k}(z)}\right\vert
 \leq \delta _{k}$ ; hence\ \par 
\quad \quad \quad $\displaystyle \ \sum_{k\in {\mathbb{N}}}{\left\vert{g_{k}(z)}\right\vert
 }\leq (1+\sum_{k\neq n}{\delta _{k}})\leq (1+\epsilon )$\ \par 
and this inequality is also satisfied for $z\in {\mathbb{T}}\backslash
 \bigcup_{n\in {\mathbb{N}}}{\Gamma _{n}}$ (since there we have
 $\ \sum_{k\in {\mathbb{N}}}{\left\vert{g_{k}(z)}\right\vert
 }\leq \epsilon $).\ \par 
\quad We now show that $T$ is bounded below:\ \par 
Let $c=(c_{n})_{n}\in \ell ^{\infty }$ and $h:=Tc=\sum_{n}^{w^{*}}{c_{n}g_{n}}$
 ; for a given $n$ we have\ \par 
\quad a) $\ \left\vert{c_{n}}\right\vert =\mathrm{m}\mathrm{a}\mathrm{x}_{z\in
 \Delta _{n}}\left\vert{c_{n}g_{n}(z)}\right\vert $\ \par 
and , for any $z\in \bar {\mathbb{D}},$ we have \ \par 
\quad b) $\ \left\vert{c_{n}g_{n}(z)}\right\vert \leq \left\vert{h(z)}\right\vert
 +\sum_{j\neq n}{\left\vert{c_{j}g_{j}(z)}\right\vert }.$\ \par 
Specializing to $z\in \Delta _{\Gamma _{n}}$ (in which case $\
 \left\vert{g_{j}(z)}\right\vert <\delta _{j}$  for $j\neq n$)\!\!\!\!
 , we obtain\ \par 
\quad \quad \quad $\displaystyle \ \left\vert{c_{n}g_{n}(z)}\right\vert \leq {\left\Vert{h}\right\Vert}+\sum_{j\neq
 n}{\delta _{j}\left\vert{c_{j}}\right\vert }\leq {\left\Vert{h}\right\Vert}+\epsilon
 {\left\Vert{c}\right\Vert}$\ \par 
and, via a), $\ \left\vert{c_{n}}\right\vert \leq {\left\Vert{h}\right\Vert}+\epsilon
 {\left\Vert{c}\right\Vert}$  for all $n,$ leading easily to
 the desired result:\ \par 
\quad \quad \quad $\displaystyle (1-\epsilon ){\left\Vert{c}\right\Vert}\leq {\left\Vert{T(c)}\right\Vert}.$\
 \par 
\quad Let us summarize what we have shown.\ \par 

\begin{Thrm}
~\label{F6}The linear mapping $T\ :\ c=(c_{n})_{n\geq N}\rightarrow
 T(c)=\sum_{n\geq N}{c_{n}g_{n}}$ implements a weak* and norm-homeomorphism
 between $\ell ^{\infty }$ and ${\mathcal{E}}=T(\ell ^{\infty
 })$ such that $\ {\left\Vert{T}\right\Vert}\leq 1+\epsilon ,\
 {\left\Vert{T^{-1}}\right\Vert}\leq (1-\epsilon )^{-1}.$
\end{Thrm}
\quad Remarks\ \par 
\quad 1. Observe that it is possible to choose $\epsilon $ so that
 the Banach-Mazur distance between $\ell ^{\infty }$ and ${\mathcal{E}}\
 ={\mathcal{E}}_{\epsilon },$ i.e., $d_{BM}(\ell ^{\infty },{\mathcal{E}}_{\epsilon
 }):=\inf \ \lbrace {\left\Vert{U}\right\Vert}{\left\Vert{U^{-1}}\right\Vert}\
 ;\ U\ isomorphism\rbrace ,$ be arbitrarily close to $1.$\ \par 
\ \par 
\quad 2. Note also that we can replace the functions $g_{n}$ by $h_{n}=u_{n}g_{n}$
  where the $u_{n}$ satisfy $\ {\left\Vert{u_{n}}\right\Vert}_{\infty
 }=1=\left\vert{u_{n}(a_{n})}\right\vert $ without affecting
 the result (that is the map $\tilde T\ :\ c=(c_{n})\rightarrow
 \sum_{n}^{w^{*}}{c_{n}h_{n}}$ still implements a weak*-homeomorphism
 between $\ell ^{\infty }$  and $\tilde T(\ell ^{\infty })$).
 Let us illustrate this in the case where the closure $F$ of
 the set $\lbrace a_{n};\ n\in {\mathbb{N}}\rbrace $ has measure
 $0$ (this will happen for instance when the sequence $(a_{n})_{n}$
  is convergent -which is in fact the situation in~\cite{BiF6}).
 This closed set $F$ is thus a peak set for $A({\mathbb{D}})$
 ; the $(u_{n})_{n}$ can then be chosen in $A({\mathbb{D}})$
 by virtue of the Rudin-Carleson interpolation theorem (cf.~\cite{BiF14}
 Chap. 2, Theorem 6.12) so as to satisfy $u_{n}(a_{k})=\delta
 _{n,k}.$ This will lead to $h_{n}(a_{k})=\delta _{n,k}$  and
 $\tilde T$ bounded below by one.\ \par 
\quad We observe that in this case the corresponding $\ell ^{1}-$ basis
 in ${\mathcal{Q}}_{{\mathcal{E}}}$ consists of evalutions $E_{a_{n}}$
 at the points $a_{n}$ (i.e., $\ {\left\langle{E_{a_{n}},h}\right\rangle}=h(a_{n}),\
 h\in H^{\infty }$) which, though not defined on all of $H^{\infty
 },$ are well-defined and weak*-continuous on $\tilde T(\ell
 ^{\infty }).$\ \par 
\ \par 
\quad 3. The construction in~\cite{BiF6} (in fact already in~\cite{BiF15})
 was inspired by~\cite{BiF16} which develops a general technique
 of producing subspaces of $H^{p}$ ($1\leq p\leq \infty $) which
 are (isomorphic to) a direct sum of infinitely many copies of
 $H^{p}.$ This technique relies on constructing subspaces ${\mathcal{M}}$
 which are $\epsilon -$ supported on a subarc $J$ of ${\mathbb{T}},$
 that is,\ \par 
\quad \quad \quad $\displaystyle \mathrm{s}\mathrm{u}\mathrm{p}_{z\in {\mathbb{T}}\backslash
 J}\left\vert{f(z)}\right\vert \leq \epsilon {\left\Vert{f}\right\Vert}_{\infty
 }$ for $f\in {\mathcal{M}}\ ;$\ \par 
it is clear that in our construction above, each (one-dimensional)
  subspace ${\mathbb{C}}g_{n}$ is $\delta _{n}-$ supported on
 the corresponding $\Gamma _{n}$ and hence our space ${\mathcal{E}}$
 above can be seen to be norm-isomorphic to $\ell ^{\infty }$
 (with a control on $d_{BM}(\ell ^{\infty },{\mathcal{E}})$ 
 similar to ours) as a consequence of~\cite{BiF16}, Lemma 1 (but
 weak*-topologies are not discussed there).\ \par 
\ \par 
\quad 4. Since for each $z\in \bar {\mathbb{D}}$ the sequence $(g_{n}(z))_{n}$
  is in $\ell ^{1},$ the expansion for $h=T(c)=\sum_{n}{c_{n}g_{n}(z)}$
 shows that $E_{z},$ the evaluation at $z,$ is well-defined on
 ${\mathcal{E}}$ and weak*-continuous, in fact $E_{z}=S^{-1}(g_{n}(z))_{n}.$
 The following continuity properties of the map $z\rightarrow
 E_{z},$ besides their intrinsic interest, will be useful later
 on. Here we denote by ${\mathcal{D}}$ the derived set of $\lbrace
 a_{n};\ n\in {\mathbb{N}}\rbrace $  (that is, the set of accumulation
 points).\ \par 

\begin{Prps}
The map $z\rightarrow E_{z}$ is continuous on $\bar {\mathbb{D}}\backslash
 {\mathcal{D}}$  and has nontangential limit at every point of
 ${\mathbb{T}}$ ; of course these limit properties transfer to
 any $h$ in ${\mathcal{E}}$ individually.
\end{Prps}
\quad Proof.\ \par 
The last statement results from the inequality\ \par 
\quad \quad \quad $\ \left\vert{h(z)-h(w)}\right\vert \leq {\left\Vert{E_{z}-E_{w}}\right\Vert}{\left\Vert{h}\right\Vert}\
 \ \ \ \ \ (h\in {\mathcal{E}},\ z,w\in \bar {\mathbb{D}}).$\ \par 
\quad Since $S$ is bicontinuous, proving the first statement amounts
 to prove the same continuity properties for the map $\Phi \
 :\ \bar {\mathbb{D}}\ni z\rightarrow S(E_{z})=(g_{n}(z))_{n}\in
 \ell ^{1}.$\ \par 
\quad {\sl Continuity of }$\Phi \ ${\sl  on }$\bar {\mathbb{D}}\backslash
 {\mathcal{D}}${\sl  }:\ \par 
Let $a\in \bar {\mathbb{D}}\backslash {\mathcal{D}}$ and $\tau
 >0$ ; by the definition of ${\mathcal{D}}$ and the fact that
 $\lbrace a_{n};\ n\in {\mathbb{N}}\rbrace $ consists of isolated
 points there exists $r>0$ such that the open disk $D_{a,r}$
 intersects at most one $\Delta _{\Gamma _{n}}.$  Consequently
 we can choose an integer $N$ such that $\displaystyle \ \sum_{n\geq
 N}{\delta _{n}}<\tau /3$  and, for $n\geq N,\ D_{a,r}\cap \Delta
 _{\Gamma _{n}}=\emptyset $ (hence $\ \left\vert{g_{n}(a)-g_{n}(z)}\right\vert
 <2\delta _{n}$  for $z\in D_{a,r}$).\ \par 
\quad Thus starting with the inequality\ \par 
\quad \quad \quad $\displaystyle \ {\left\Vert{\Phi (a)-\Phi (z)}\right\Vert}\leq
 \sum_{n<N}{\left\vert{g_{n}(a)-g_{n}(z)}\right\vert }+\sum_{n\geq
 N}{\left\vert{g_{n}(a)-g_{n}(z)}\right\vert }$\ \par 
we obtain, for $z\in D_{a,r}$\ \par 
\quad \quad \quad $\displaystyle \ {\left\Vert{\Phi (a)-\Phi (z)}\right\Vert}\leq
 \sum_{n<N}{\left\vert{g_{n}(a)-g_{n}(z)}\right\vert }+2\tau /3.$\ \par 
\quad Now the continuity (with respect to $z$) of the first term on
 the right hand side of this inequality enables us to pick $\eta
 \ (0<\eta <r)$ such that, for $\ \left\vert{z-a}\right\vert
 <\eta ,\ \sum_{n<N}{\left\vert{g_{n}(a)-g_{n}(z)}\right\vert
 }<\tau /3,$ yielding for such $z,$ $\ {\left\Vert{\Phi (a)-\Phi
 (z)}\right\Vert}<\tau .$ \ \par 
\quad {\sl Nontangential continuity of }$\Phi ${\sl  : }\ \par 
Let $a\in {\mathbb{T}}$ and, for $0<r<1,$ let $St_{r}$ the Stoltz
 domain defined by $a$ and $r$ (that is $St_{r}$ is the interior
 of the convex hull of $r{\mathbb{D}}\cup \lbrace a\rbrace $
 ) ; there is at most one $n$ (say $n_{0}$) such that $a\in \Delta
 _{r_{n}}$ ; elementary geometric considerations show that, as
 soon as the length of $\Gamma _{n}$ is smaller than $2{\sqrt{1-r^{2}}}$
 (and $n\neq n_{0}$) we have $St_{r}\cap \Delta _{\Gamma _{n}}=\emptyset
 .$ This done we proceed as before to conclude. $\displaystyle
 \hfill\blacksquare $\ \par 
\ \par 
5. We note also that this construction provides (by restriction
 of the map $T$ to $c_{0}$) a norm-isomorphism between $c_{0}$
 and a subspace of the disk algebra, in other words a subspace
 of $A({\mathbb{D}})$  with a $c_{0}-$ basis.\ \par 
\ \par 
6. Finally in the case when $(a_{n})_{n}$ is convergent the set
 ${\mathcal{D}}$ is reduced to a singleton, say $\lbrace a\rbrace
 ,$ we observe that the subspace ${\mathcal{E}}$ does not seem
 to be "missing by much" its inclusion in the disk algebra $A({\mathbb{D}})$
 : indeed only continuity at the point $a$ might be lacking.
 Nevertheless this is enough to substantially "enlarge" the space
 ${\mathcal{E}}$ since (being isomorphic to $\ell ^{\infty }$)
 it is not even separable.\ \par 

\subsection{A construction with interpolating sequences.}
\ \par 
\quad As is well-known, a sequence $(b_{n})_{n}$ in ${\mathbb{D}}$
 is said to be an interpolating sequence if the mapping ${\mathcal{J}}\
 :\ h\in H^{\infty }\rightarrow (h(b_{n}))_{n}\in \ell ^{\infty
 }$ is onto.\ \par 
\quad Here the construction goes even faster thanks to the following
 result of Pehr Beurling.\ \par 

\begin{Thrm}
~\cite{BiF17}  ~\label{F11}Given an interpolating sequence, $(b_{n})_{n}$
 in ${\mathbb{D}},$ there exists a sequence of bounded analytic
 functions $(\beta _{n})_{n\in {\mathbb{N}}}$  having the following
 properties:\par 
\quad (i) $\beta _{k}(b_{j})=\delta _{k,j}$  for all $j,k\in {\mathbb{N}}$ and \par 
\quad (ii) $\mathrm{s}\mathrm{u}\mathrm{p}_{z\in {\mathbb{D}}}\sum_{n\in
 {\mathbb{N}}}{\left\vert{\beta _{n}(z)}\right\vert }<M$  (where
 $M$ is the interpolating constant).
\end{Thrm}
\quad Indeed, since Hypothesis (ii) is exactly Assertion (ii) of Proposition~\ref{FB0}
 we have already a weak*-continuous linear map $T$ from $\ell
 ^{\infty }$ in $H^{\infty }$ such that $\forall n\in {\mathbb{N}},\
 T\epsilon _{n}=\beta _{n}.$ Moreover, for any $c\in \ell ^{\infty
 }$ and any $n\in {\mathbb{N}},$ we have\ \par 
\quad \quad \quad $\displaystyle \ \left\vert{c_{n}}\right\vert =\left\vert{c_{n}\beta
 _{n}(a_{n})}\right\vert =\left\vert{\sum_{k}{c_{k}\beta _{k}(a_{n})}}\right\vert
 \leq {\left\Vert{T(c)}\right\Vert}_{\infty }$\ \par 
which shows that the map $T$ is bounded below by $1.$ Thus, by
 Prop.~\ref{F10}, ${\mathcal{E}}:=T(\ell ^{\infty })$ is weak*-closed
 and $T,$ seen as a map from $\ell ^{\infty }$ into ${\mathcal{E}},$
 is a weak*-homeomorphism whose inverse is obviously the restriction
 to the subspace ${\mathcal{E}}$ of the map ${\mathcal{J}}$ (defined
 at the beginning of this subsection). Moreover by (i) we can
 write for any $a=(a_{n})_{n}\in \ell ^{1}$ and $c=(c_{n})_{n}\in
 \ell ^{\infty }$\ \par 
\quad \quad \quad $\displaystyle \ {\left\langle{a,c}\right\rangle}=\sum_{n}{a_{n}c_{n}}=\sum_{j,k}{a_{j}c_{k}\delta
 _{j,k}}={\left\langle{\sum_{j}{a_{j}E_{b_{j}},T(c)}}\right\rangle}.$\ \par 
\quad Thus the "preadjoint" $S\ (\in {\mathcal{L}}({\mathcal{Q}}_{{\mathcal{E}}},\
 l^{1}))$  of $T$ is defined by $S(\sum_{j}{a_{j}E_{b_{j}}})=a$
 : in other words, the predual ${\mathcal{Q}}_{{\mathcal{E}}}$
 admits $(E_{b_{n}})_{n}$ as an $\ell ^{1}-$ basis and ${\mathcal{E}}$
 admits $(\beta _{n})_{n}$  as an $\ell ^{\infty }-$ basis.\ \par 
\quad Summing up we have shown:\ \par 

\begin{Thrm}
~\label{F14}Let $(b_{n})_{n}$ be an interpolating sequence in
 ${\mathbb{D}}$ and $(\beta _{n})_{n}$ a sequence of $H^{\infty
 }$ functions satisfying (i) and (ii) of Theorem~\ref{F11} ;
 then the linear mapping\par 
\quad \quad \quad $T\ :\ c=(c_{n})_{n\in {\mathbb{N}}}\rightarrow T(c)=\sum_{n\in
 {\mathbb{N}}}{c_{n}\beta _{n}}$ \par 
implements a weak* and norm-homeomorphism between $\ell ^{\infty
 }$  and ${\mathcal{E}}=T(\ell ^{\infty })$ (with inverse ${\mathcal{J}}_{\mid
 {\mathcal{E}}}$) such that $\ {\left\Vert{T}\right\Vert}\leq
 M,\ {\left\Vert{T^{-1}}\right\Vert}\leq 1$ ; the sequence $(E_{b_{n}})_{n}$
  is an $\ell ^{1}-$ basis for ${\mathcal{Q}}_{{\mathcal{E}}}$
 and $(\beta _{n})_{n}$ an $\ell ^{\infty }-$ basis for ${\mathcal{E}}.$
\end{Thrm}
\quad It turns out that the reference~\cite{BiF17} is a preprint which
 is far from widely available. Thus we give below a proof of
 a weaker version (bound $M^{2}$ instead of $M$) based on a lemma
 of Drury~\cite{BiF18} and a more classical formulation of results
 on interpolating sequences.\ \par 
\quad We start with a "finite" version of what is needed. This is in~\cite{BiF19},
 Prop.~\ref{F10}, (with $n=1$) where Drury's lemma is used along
 the lines of~\cite{BiF20}. We include the proof for the sake
 of completeness.\ \par 

\begin{Lmm}
~\label{F15}Suppose the sequence $(b_{j})_{1\leq j\leq n}$ consisting
 of pairwise distinct elements of ${\mathbb{D}}$ is interpolating
 with constant $M$ ; then we can find bounded analytic fonctions
 $(\beta _{j})_{1\leq j\leq n}$  such that \par 
\quad (i) $\beta _{k}(b_{j})=\delta _{k,j}$ $1\leq k\neq j\leq n$ and\par 
\quad (ii) $\mathrm{s}\mathrm{u}\mathrm{p}_{z\in {\mathbb{D}}}\sum_{1\leq
 j\leq n}{\left\vert{\beta _{j}(z)}\right\vert }<M^{2}.$
\end{Lmm}
\quad Proof.\ \par 
Let $\lambda $ be a primitive $n-$ th roots of unity in ${\mathbb{C}}$
 ; the hypothesis ($(b_{j})_{1\leq j\leq n}$  interpolating with
 constant $M$) ensures the existence of functions $(v_{j})_{1\leq
 j\leq n}$  such that $v_{k}(b_{j})=\lambda ^{kj}$ $1\leq j\leq
 n$ and $\mathrm{m}\mathrm{a}\mathrm{x}_{1\leq j\leq n}{\left\Vert{v_{j}}\right\Vert}_{\infty
 }\leq M.$\ \par 
\quad Define the functions $\varphi _{j}$ on ${\mathbb{D}}$ by\ \par 
\quad \quad \quad $\displaystyle \varphi _{j}(z)=\frac{1}{n}\sum_{k=1}^{n}{\lambda
 ^{-kj}v_{k}(z)},\ \ \ 1\leq j\leq n,$\ \par 
the "Fourier transform" of the function $\displaystyle v(\cdot
 ,z):=v_{k}(z),$ where $z\in {\mathbb{D}}$ is a parameter, on
 the group $G$ of the $n-$ th root of unity with the measure
 $\displaystyle \ \frac{1}{n}\sum_{k=1}^{n}{\delta }_{k}.$\ \par 
Observe that\ \par 
\quad \quad \quad $\displaystyle \varphi _{j}(b_{l})=\frac{1}{n}\sum_{k=1}^{n}{\lambda
 ^{-kj}v_{k}(b_{l})}=\frac{1}{n}\sum_{k=1}^{n}{\lambda ^{(l-j)k}=\delta
 _{j,l}},\ 1\leq j\neq l\leq n.$\ \par 
\quad Clearly, the $\varphi _{j}$ 's belong to $H^{\infty }$ (and $\
 {\left\Vert{\varphi _{j}}\right\Vert}_{\infty }\leq M$) ;  a
 straigthforward computation shows that, for any $z$ in ${\mathbb{D}},$\ \par 
\quad \quad \quad $\displaystyle \ \left\vert{\varphi _{j}(z)}\right\vert ^{2}=\frac{1}{n^{2}}(\sum_{k=1}^{n}{\left\vert{\varphi
 _{k}(z)}\right\vert ^{2}}+\sum_{1\leq k\neq l\leq n}{\lambda
 ^{j(l-k)}A_{k,l}})$\ \par 
with with $A_{k,l}=v_{k}(z){\overline{v_{l}(z)}}.$\ \par 
\quad Since, for given $l\neq k$ (as already observed before), it follows that\ \par 
\quad \quad \quad $\displaystyle \ \sum_{k=1}^{n}{\left\vert{\varphi _{k}(z)}\right\vert
 ^{2}}=\frac{1}{n^{2}}\sum_{j=1}^{n}{\sum_{k=1}^{n}{\left\vert{v_{k}(z)}\right\vert
 ^{2}}}\leq M^{2}.$\ \par 
\quad The above equality can  also be easily seen as a consequence
 of Plancherel theorem between the group on $n$-th roots of unity
 in ${\mathbb{C}}$ and its dual.\ \par 
\quad The proof of the lemma is thus completed by setting $\beta _{j}=(\varphi
 _{j})^{2}.$ $\displaystyle \hfill\blacksquare $\ \par 
\ \par 
\quad Thus, given an interpolating sequence $(b_{n})_{n\in {\mathbb{N}}}$
 with constant $M,$ for each $n$ we have a finite sequence of
 analytic functions $\displaystyle (\beta _{j,n})_{1\leq j\leq
 n}$ satisfying (i) and (ii) of lemma~\ref{F15}. Now a standard
 procedure completes the proof of (the weak version of) P. Beurling's
 theorem.\ \par 
\quad Since this procedure requires successive extractions of subsequences
 it is convenient to introduce the following notation for (strictly
 increasing in our case) maps $\varphi $ and $\psi $ of ${\mathbb{N}}$
 in itself:\ \par 
\quad $\varphi \prec \psi $ if there exists another strictly increasing
 sequence $j\ :\ {\mathbb{N}}\rightarrow {\mathbb{N}}$ such that
 $\psi =\varphi \circ j$ (thus for any sequence whatsoever $(u_{n})_{n}$
 the sequence $(u_{\psi (n)})_{n}$  is a subsequence of the sequence
 $(u_{\varphi (n)})_{n}$ (itself subsequence of $(u_{n})_{n}$ ).\ \par 
\quad Note that all the $\beta _{j,n}$ are of norm bounded by $M^{2}.$
 As a first step we extract from the bounded sequence $(\beta
 _{1,n})_{n}$  a subsequence $(\beta _{1,\varphi _{1}(n)})_{n}$
 which is weak*-convergent to say $\beta _{1}.$ Without loss
 of generality we may and do assume that $\varphi _{1}(1)\geq
 2$ so that the subsequence $(\beta _{2,\varphi _{1}(n)})_{n}$
  of $(\beta _{2,n})_{n}$ is well-defined. As a second step we
 extract from $(\beta _{2,\varphi _{1}(n)})_{n}$ a subsequence
 $(\beta _{2,\varphi _{2}(n)})_{n}$  (that is, $\varphi _{1}\prec
 \varphi _{2},$ with again, as precaution for the next step,
 $\varphi _{2}(1)\geq 3$ ) which is weak*-convergent to say $\beta
 _{2}.$ Continuing in this fashion we clearly obtain a sequence
 of strictly increasing maps $\varphi _{k}$ from ${\mathbb{N}}$
  to ${\mathbb{N}}$ such that $\varphi _{1}\prec \varphi _{2}\prec
 \cdot \cdot \cdot \prec \varphi _{k}\prec \cdot \cdot \cdot
 $ with, for all $k,\ \varphi _{k}(1)\geq k+1$ and the sequence
 $(\beta _{k,\varphi _{k}(n)})_{n}$ is weak*-convergent to say
 $\beta _{k}.$\ \par 
\quad We now show that the sequence $(\beta _{j})_{j}$ satisfies the
 desired properties:\ \par 
let $j,l$ in ${\mathbb{N}}$ ; by the initial definition of the
 (finite)  sequence$\lbrace \beta _{j,\varphi _{j}(n)}\rbrace
 _{1\leq j\leq \varphi _{j}(n)},$ we have $\beta _{j,\varphi
 _{j}(n)}(b_{l})=\delta _{j,l},$ an equality which is of course
 preserved by taking the limit as $n\rightarrow \infty ,$ yielding
 thus $\beta _{j}(b_{l})=\delta _{j,l}.$\ \par 
\quad Similarly we have for any given $z\in {\mathbb{D}},\ \sum_{1\leq
 j\leq \varphi _{j}(n)}{\left\vert{\beta _{j,\varphi _{j}(n)}}\right\vert
 }<M^{2}\ $ ; in particular for any integer $N$ we will have,\ \par 
\quad \quad \quad $\displaystyle \ \sum_{1\leq j\leq N}{\left\vert{\beta _{j,\varphi
 _{j}(n)}}\right\vert }<M^{2}$ as soon as $\varphi _{j}(n)\geq N,$\ \par 
and taking the limit as $\displaystyle n\rightarrow \infty $
 we get $\ \sum_{1\leq j\leq N}{\left\vert{\beta _{j}}\right\vert
 }<M^{2}$ from which the proof is easily concluded. $\displaystyle
 \hfill\blacksquare $\ \par 
\ \par 

\begin{Rmrq}
Of course, the above construction is more "existential" than
 "explicit". However, as indicated in the introduction, in some
 situations it is possible to obtain a "concrete" construction
 (see Theorem~\ref{F8}).
\end{Rmrq}
\ \par 

\subsection{Additional remarks.}
\ \par 
\quad When the elementary construction (developed in section~\ref{F12})
 can be modified via the Rudin-Carleson interpolation theorem
 we obtain a weak*-copy ${\mathcal{E}}$ of $\ell ^{\infty }$
 whose predual ${\mathcal{Q}}_{{\mathcal{E}}}$  admits an $\ell
 ^{1}-$ basis consisting of point evaluations at the points $a_{n}$
 (cf. Remark 2 following Theorem~\ref{F6} above). It is interesting
 to observe that in such a case the predual ${\mathcal{Q}}_{{\mathcal{E}}}$
 admits also $\ell ^{1}-$ bases of the form $(E_{b_{n}})_{n}$
 with the $b_{n}$'s in ${\mathbb{D}}.$ Thus ${\mathcal{E}}$ admits
 the system of functions biorthogonal to $(E_{b_{n}})_{n},$ say
 $(\beta _{n})_{n},$  as an $\ell ^{\infty }-$ basis and it follows
 easily that $(b_{n})_{n}$ is an interpolating sequence. Before
 stating a precise result we recall a Paley-Wiener type stability
 theorem that will be the key to its proof.\ \par 

\begin{Thrm}
~\label{F4} (~\cite{BiF13} Theorem 9.1, p. 84) Let $(x_{n})_{n}$
 be a Schauder basis in a Banach space $X.$\par 
\quad a) A sequence $(y_{n})_{n}$ in $X$ for which there exists a constant
 $\lambda \in \rbrack 0,1\lbrack $ such that\par 
\quad \quad \quad $\displaystyle \forall \alpha =(\alpha _{j})_{j}\in c_{00},\
 {\left\Vert{\sum_{j}{\alpha _{j}(x_{j}-y_{j})}}\right\Vert}\leq
 \lambda {\left\Vert{\sum_{j}{\alpha _{j}x}_{j}}\right\Vert}$\par 
is a basis (equivalent to the basis $(x_{n})_{n}$).\par 
\quad b) If in addition $(x_{n})_{n}$ is an $\ell ^{1}-$ basis then
 there exists $\eta >0$ such that $(y_{n})_{n}$ is an $\ell ^{1}-$
 basis whenever $\ {\left\Vert{x_{n}-y_{n}}\right\Vert}<\eta $ for all $n.$
\end{Thrm}
\quad Proof.\ \par 
\quad a) The hypothesis implies that the linear map $\ \sum_{j}{\alpha
 _{j}x_{j}}\rightarrow \sum_{j}{\alpha _{j}(x_{j}-y_{j})}$ defined
 on the linear span of the $x_{n}$'s extends into a bounded linear
 operator $U\in {\mathcal{L}}(X)$ satisfying $\ {\left\Vert{U}\right\Vert}\leq
 \lambda <1.$ Consequently the operator $A:=I_{X}-U$ is invertible
 and since $Ax_{n}=y_{n}$  for all $n$ we have the desired conclusion.\ \par 
\quad b) Since $(x_{n})_{n}$ is an $\ell ^{1}-$ basis there exists
 a constant $\nu $ such that for all $\alpha =(\alpha _{j})_{j}\in
 c_{00}$ we have $\ {\left\Vert{\alpha }\right\Vert}_{1}\leq
 \nu {\left\Vert{\sum_{j}{\alpha _{j}x}_{j}}\right\Vert}$ ; for
 $\alpha \in c_{00}$ and $(y_{n})_{n}$ such that $\ {\left\Vert{x_{j}-y_{j}}\right\Vert}\leq
 \eta $  for all $j,$ the inequalities\ \par 
\quad \quad \quad $\displaystyle \ {\left\Vert{\sum_{j}{\alpha _{j}(x_{j}-y_{j})}}\right\Vert}\leq
 \eta {\left\Vert{\alpha }\right\Vert}_{1}\leq \eta \nu {\left\Vert{\sum_{j}{\alpha
 _{j}x}_{j}}\right\Vert}$\ \par 
enable us to apply a) whenever $\eta \nu <1.$ $\displaystyle
 \hfill\blacksquare $\ \par 
\ \par 

\begin{Thrm}
~\label{F7}Let ${\mathcal{E}}$ be a weak*-closed subspace of
 $H^{\infty }$ consisting of functions continuously extendable
 to a subset $\Omega $ of $\bar {\mathbb{D}}\ (\Omega \supset
 {\mathbb{D}})$ such that the map $\Omega \ni z\rightarrow E_{z}$
  is valued in ${\mathcal{Q}}_{{\mathcal{E}}}$ and continuous.
 Suppose furthermore that ${\mathcal{Q}}_{{\mathcal{E}}}$ admits
 an $\ell ^{1}-$ basis of the form $(E_{a_{n}})_{n}$ with the
 $a_{n}$'s in ${\mathbb{T}}\cap \Omega .$ Then there exists (strictly)
 positive numbers $\tau _{n}$ such that\par 
\quad (i) the closed disks $\bar {\mathbb{D}}_{n}:=\bar {\mathbb{D}}_{a_{n},\tau
 _{n}}$ are pairwise disjoints, and \par 
\quad (ii) for any sequence $(b_{n})_{n}$ chosen so that (for all $n$)
  $b_{n}\in {\mathbb{D}}\cap {\mathbb{D}}_{n}$ the sequence $(E_{b_{n}})_{n}$
  is an $\ell ^{1}-$ basis for ${\mathcal{Q}}_{{\mathcal{E}}}$
 (and consequently $(b_{n})_{n}$ is an interpolating sequence).
\end{Thrm}
(Of course in the context of this theorem the system $(g_{n})_{n}$
 in ${\mathcal{E}},$ biorthogonal to $\displaystyle (E_{a_{n}})_{n}$
 is an $\ell ^{\infty }$-basis for ${\mathcal{E}}$ and the map
 implementing the weak*-homeomorphism between $\displaystyle
 \ell ^{\infty }$ and ${\mathcal{E}}$ is given by $T(\gamma )=\sum_{j}{\gamma
 _{j}g_{j}}.$)\ \par 
\quad Proof.\ \par 
\quad (i) Since $(E_{a_{n}})_{n}$ is an $\ell ^{1}-$ basis in ${\mathcal{Q}}_{{\mathcal{E}}},$
 we have $\mathrm{i}\mathrm{n}\mathrm{f}_{n\neq m}{\left\Vert{E_{a_{n}}-E_{a_{m}}}\right\Vert}>0$
 ; this inequality combined with the continuity of the map $z\rightarrow
 E_{z}$  at the points $a_{k}$ ensures that for all $n,\ d(a_{n},\lbrace
 a_{m}\ ;\ m\neq n\rbrace >0.$\ \par 
\quad (ii) : Immediate consequence of the continuity of the map $z\rightarrow
 E_{z}$  combined with the above Theorem~\ref{F4} (with $X={\mathcal{Q}}_{{\mathcal{E}}}$
 , $(x_{n})_{n}=(E_{a_{n}})_{n},\ $ $(y_{n})_{n}=(E_{b_{n}})_{n}$).
 $\displaystyle \hfill\blacksquare $\ \par 
\ \par 

\begin{Rmrq}
This last argument is in fact a particular case of the Krein-Milman-Rutman
 stability theorem (cf.~\cite{BiF13}  Theorem 10.3 p. 98) simplified
 here because we deal with the $\ell ^{1}-$ norm.
\end{Rmrq}
\ \par 
\quad As announced in the introduction this rich supply of interpolating
 sequences in the vicinity of $(a_{n})_{n}$ can be obtained under
 the mere hypothesis of separation of the $a_{n}$'s on the unit
 circle. This is probably known to specialists in this area but
 for the sake of completeness we give a precise statement and
 prove it using the Carleson's characterization of interpolating
 sequences, namely a sequence $(b_{n})_{n}$ in ${\mathbb{D}}$
 is interpolating if there exists $m>0$ such that\ \par 
\quad \quad \quad $\displaystyle \inf \ _{k\in {\mathbb{N}}}\prod_{j\neq k}{\left\vert{\frac{b_{k}-b_{j}}{1-\bar
 b_{k}b_{j}}}\right\vert }>m.$\ \par 
\quad We will even exhibit a concrete linear extension operator associated
 to any such sequence $(b_{n})_{n}$ close to $(a_{n})_{n}$ (but
 of course the weak*-closed subspace of $H^{\infty }$ -image
 of this operator- depends on the sequence $(b_{n})_{n}$ while
 in Theorem~\ref{F7} we obtain the same subspace). We will use
 the standard terminology $B,\ B_{k},\ k=1,2,...$ denote the
 usual Blaschke products,\ \par 
\quad \quad \quad $\displaystyle B(z):=\prod_{n}{\frac{b_{n}-z}{1-\bar b_{n}z}},\
 B_{k}(z):=\prod_{n\neq k}{\frac{b_{n}-z}{1-\bar b_{n}z}},\ k\in
 {\mathbb{N}}$\ \par 
and, following~\cite{BiF8}  formula 7.6, p.48, we introduce the
 functions $g_{k}$ :\ \par 
\quad \quad \quad $\displaystyle g_{k}(z)=\frac{1-\left\vert{b_{k}}\right\vert
 ^{2}}{1-\bar b_{k}z}\frac{B_{k}(z)}{B_{k}(b_{k})},\ k\in {\mathbb{N}}.$\ \par 
\ \par 

\begin{Thrm}
~\label{F8}Let $(a_{n})_{n}$ be a separated sequence in ${\mathbb{T}}$
 and let $m\in \rbrack 0,1\lbrack $ ; then, there exists a sequence
 $(\eta _{n})_{n}$ of strictly positive numbers such that any
 sequence $(b_{n})_{n}$ in ${\mathbb{D}}$  satisfying $\ \left\vert{b_{n}-a_{n}}\right\vert
 <\eta _{n}$ for all $n$ is an interpolating sequence with constant
 no smaller than $m.$ Moreover the functions $(g_{n})_{n}$ associated
 to the sequence $(b_{n})_{n}$ by the above formula generate
 a weak*-copy of $\ell ^{\infty }$ (with $(g_{n})_{n}$  as an
 $\ell ^{\infty }-$ basis).
\end{Thrm}
\quad Proof.\ \par 
We choose first positive numbers $d_{n}$ such that the closed
 disks $\bar {\mathbb{D}}_{a_{n},d_{n}}$ are pairwise disjoint.\ \par 
{\sl Notational simplification for this proof:} $\displaystyle
 \bar {\mathbb{D}}_{a,d}$ means $\bar {\mathbb{D}}\bigcap \bar
 {\mathbb{D}}_{a,d}.$\ \par 
\quad We get by induction on $N$ a sequence of positive numbers $\tau
 _{n}$ and a sequence of compact sets $F_{n}$ in $\bar {\mathbb{D}}^{n}$
 in the following fashion.\ \par 
\quad Set $\tau _{1}=d_{1}$ and $F_{1}=\bar {\mathbb{D}}_{a_{1},\tau
 _{1}}$ ; the map\ \par 
\quad \quad \quad $\displaystyle \Theta \ :\ \bar {\mathbb{D}}_{a_{2},d_{2}}{\times}F_{1}\ni
 (t,s)\rightarrow \left\vert{\frac{t-s}{1-\bar ts}}\right\vert $\ \par 
is continuous and satisfies $\Theta (a_{2},s)\ (=1)\ >m$ for
 all $s.$ Thus the open set $\Theta ^{-1}(\rbrack m,\infty \lbrack
 )$ contains the (compact) set $\lbrace a_{2}\rbrace {\times}F_{1}$
 ; therefore, by a standard (and elementary)  compacity argument,
 it contains a (compact) set $F_{2}$ of the form $F_{2}=\bar
 {\mathbb{D}}_{a_{2},\tau _{2}}$ with $0<\tau _{2}(\leq d_{2}).$\ \par 
\quad Suppose that $\tau _{1},\tau _{2},..,\tau _{N}$ have been found
 such that for any $s=(s_{1},..,s_{N})\in F_{N}:=\prod_{1\leq
 n\leq N}{\bar {\mathbb{D}}_{a_{n},d_{n}}},$  we have\ \par 
\quad \quad \quad $\displaystyle \Theta _{N}(s):=\prod_{j\neq k,\ j,k\leq N}{\left\vert{\frac{s_{j}-s_{k}}{1-\bar
 s_{j}s_{k}}}\right\vert }>m.$\ \par 
\quad The function\ \par 
\quad \quad \quad $\displaystyle \Theta _{N+1}\ :\ \bar {\mathbb{D}}_{a_{N+1},d_{N+1}}{\times}F_{N}\ni
 (t,s)\rightarrow (\prod_{1\leq j\leq N}{\left\vert{\frac{t-s_{j}}{1-\bar
 ts_{j}}}\right\vert }){\times}\Theta _{N}(s)),$\ \par 
is continuous and satisfies $\Theta _{N+1}(a_{N+1},s)=\Theta
 _{N}(s)>m$  for all $s$ in $F_{N}.$ The same continuity-compacity
 argument as above yields the existence of a (compact) set $F_{N+1}=\bar
 {\mathbb{D}}_{a_{N+1},\tau _{N+1}}{\times}F_{N}$ such that,
 for any $(t,s)\in F_{N+1},\ \Theta _{N+1}(t,s)\ (=1)\ >m.$ Thus
 the first statement of our theorem holds for any choice of the
 $\eta _{n}$' s satisfying $\eta _{n}\leq \tau _{n}.$\ \par 
\quad To prove the last assertion we just have to show that the functions
 $(g_{k})_{k}$ satisfy Condition (ii) of Theorem~\ref{F11} (indeed,
 as well as the functions $(B_{k}(z_{k}(b_{k}))_{k},$ they satisfy
 Condition (i) which in turn ensures that the associated map
 $(c_{n})_{n}\rightarrow \sum_{n}{c_{n}g_{n}}$ is bounded below).
 By the Carleson condition we have $\ \left\vert{B_{k}(b_{k})}\right\vert
 \geq m$  (we assume of course that all sequences $(b_{n})_{n}$
 considered here satisfy $\ \left\vert{b_{n}-a_{n}}\right\vert
 <\tau _{n}$) hence $\ \left\vert{B_{k}(z_{k}(b_{k}))}\right\vert
 \leq 1/m$  for all $k$ in ${\mathbb{N}}$ and $z$ in ${\mathbb{D}}.$\ \par 
\quad Therefore, we just need to perform some work on the additional
 factor  $\ \frac{1-\left\vert{b_{n}}\right\vert ^{2}}{\left\vert{1-\bar
 b_{n}z}\right\vert }$ to ensure (ii). Firstly we note that for
 any $z,t$ in ${\mathbb{D}}$ we have by geometrical observation
 the inequality $\ \left\vert{1-\bar b_{n}z}\right\vert \geq
 1-\left\vert{b_{n}}\right\vert $  which leads to $\ \frac{1-\left\vert{b_{n}}\right\vert
 ^{2}}{\left\vert{1-\bar b_{n}z}\right\vert }\leq 2.$ Next we
 let $(\epsilon _{n})_{n}$ another sequence of positive numbers
 such that $\Sigma _{j}\epsilon _{j}=\epsilon <m.$ We apply again
 the standard continuity-compacity argument to the function $\varphi
 :\ K_{n}{\times}{\mathbb{D}}_{a_{n},\tau _{n}}\ni (z,t)\rightarrow
 \frac{1-\left\vert{t}\right\vert ^{2}}{\left\vert{1-\bar tz}\right\vert
 }$  (obviously continuous and with value $0$ on $K_{n}{\times}\lbrace
 a_{n}\rbrace $)  to obtain $\eta _{n}>0$ (and of course $\eta
 _{n}<\tau _{n}$) such that $\varphi $ is smaller than $\epsilon
 _{n}$ on $K_{n}{\times}{\mathbb{D}}_{a_{n},\eta _{n}}.$\ \par 
\quad Now let $z\in {\mathbb{D}}$ :\ \par 
\quad a) if $z$ does not belong to the union of the disks ${\mathbb{D}}_{a_{n},\tau
 _{n}}$  then, for each $n,\ \left\vert{g_{n}(z)}\right\vert
 \leq \epsilon _{n}/m$ (recall that $\ \left\vert{B_{n}(b_{n})}\right\vert
 >m$) and hence $\ \sum_{j}{\left\vert{g_{j}(z)}\right\vert }<\epsilon
 $  ; \ \par 
\quad b) otherwise $z$ belongs to exactly one of these disks say ${\mathbb{D}}_{a_{k},\tau
 _{k}}$ ; in this case we get\ \par 
\quad \quad \quad $\displaystyle \ \sum_{j}{\left\vert{g_{j}(z)}\right\vert }=\sum_{j\neq
 k}{\left\vert{g_{j}(z)}\right\vert }+\left\vert{g_{k}(z)}\right\vert
 \leq \sum_{j\neq k}{\epsilon _{j}}+2/m.$\ \par 
\quad Thus Condition (ii) is clearly satisfied ; this concludes the
 proof. $\displaystyle \hfill\blacksquare $\ \par 
\ \par 
\quad In the context of Theorem~\ref{F7} one might think that the closure
 of the set $\lbrace a_{n};\ n\in {\mathbb{N}}\rbrace $ is necessarily
 of measure $0.$ In fact it is easy to sketch examples showing
 that this idea is wrong.\ \par 
\ \par 

\begin{Exmp}
We start with a closed interval $J=\lbrack \alpha ,\beta \rbrack
 $ of the real line and build in this interval a Cantor type
 set $K=\bigcap _{k\geq 1}J_{k}$  of positive measure in the
 standard fashion, that is, $J_{k}=J_{k-1}\backslash \bigcup_{n=2^{k-1}}^{2^{k}-1}{I_{n}}$
 (with $J_{0}=J$) where, for $2^{k-\mathrm{ }\mathrm{1}}\leq
 n<2^{k},\ I_{n}=\rbrack b_{n},c_{n}\lbrack $  is an open interval
 of length $\rho _{k}$ (with $\ \sum_{k>0}{\rho _{k}}<\beta -\alpha
 $  to ensure that $K$ has positive measure) removed in the usual
 "centered" way from each of the successive $2^{k-1}$ intervals
 whose union is $J_{k-1}.$ It follows easily from this classical
 construction that for any sequence $(x_{n})_{n}$ satisfying
 $b_{n}\leq x_{n}\leq c_{n}\ (n\geq 1)$  the endpoints $b_{n},\
 c_{n}$ are in the closure of the set $\lbrace x_{j};\ j\geq
 1\rbrace \ ;$  thus this closure contains the set $K$ (in fact
 it is equal to $K\bigcup \lbrace x_{j};\ j\geq 1\rbrace $).
 Now choosing $x_{n}$ in the open interval $\rbrack b_{n},c_{n}\lbrack
 $ we have $d(x_{n},\lbrace x_{j};\ j\neq n\rbrace )\ (=\mathrm{m}\mathrm{i}\mathrm{n}\lbrace
 x_{n}-b_{n},c_{n}-x_{n}\rbrace \ )\ >0.$  Taking $\beta -\alpha
 <2\pi $ the sequence $(a_{n})_{n}:=(e^{ix_{n}})_{n}$ will be
 one of the desired examples.
\end{Exmp}
\ \par 

\subsection{Weak*-copies of $\ell ^{\infty }$ are w*-complemented.~\label{F13}}
\ \par 
\quad As shown in the theorem~\ref{F14}, in the second type of construction
 the map $T^{-1}(:\ {\mathcal{E}}\rightarrow \ell ^{\infty })$
  is exactly the restriction to ${\mathcal{E}}$ of the interpolating
 map ${\mathcal{J}}\ :\ H^{\infty }\rightarrow \ell ^{\infty
 }$ ${\mathcal{J}}(h):=(h(b_{n}))_{n}.$ Observe also that the
 map $P\ :\ H^{\infty }\rightarrow H^{\infty }$  defined by $P(h):=T({\mathcal{J}}(h))$
 is a weak*-continuous projection (norm-bounded by $M^{2}$) whose
 range is the space ${\mathcal{E}}$  and, consequently, ${\mathcal{E}}$
 is weak*-complemented.\ \par 
This fact is no "accident" and is valid in a more general setting
 which we now describe. Let $X$ denote an arbitrary Banach space
 $X$ and let $Z=X^{*}.$ Note first that the map ${\mathcal{J}}$
  above (where $X=L^{1}/H_{0}^{1}$)  can be defined by ${\mathcal{J}}(h)=({\left\langle{E_{b_{j}},h}\right\rangle})_{j}$
  (where $E_{b_{j}}$ is the element in $L^{1}/H_{0}^{1}$  evaluation
 at $b_{j}$)\!\!\!\! .\ \par 
\quad We generalize this map by associating to any (bounded) sequence
 $K=(k_{j})_{j}$ of elements in $X$ an "interpolation map" $R\
 (=R_{K})$ from $Z$ in $\ell ^{\infty }$ defined by $R(h)=({\left\langle{k_{j},h}\right\rangle})_{j}\
 ,\ \ \ h\in Z.$\ \par 
\quad We also associate to the sequence $K$ the (bounded linear)  map
 $V\ (=V_{K})$ from $\ell ^{1}$ into $X$ defined by $V(\gamma
 )=\sum_{n}{\gamma _{n}k_{n}},\ \ \gamma \in \ell ^{1}$ and note
 that $V^{*}=R.$\ \par 
\quad We will say that $K$ is an $\ell ^{\infty }-$ interpolating sequence
 ($\displaystyle \ell ^{\infty }$-IS for short) if the map $R_{K}$
 is onto (or equivalently if the map $V_{K}$  is bounded below).We
 will denote by $\displaystyle X_{K}$ the closure of the range
 of $\displaystyle V_{K}.$\ \par 
\quad In this setting we have the following general result which is,
 in fact, with a slightly different formulation, the equivalence
 of Assertions 2 and 3 of Theorem 10 in~\cite{BiF23}, page 48,
 (see below, Remark 3).\ \par 

\begin{Thrm}
~\label{FB2}a) Let ${\mathcal{E}}$ be a (linear) subspace of
 $Z$ weak*-homeomorphic to $\ell ^{\infty }$ (via $T$). Then\par 
\quad (i) there exists a (bounded) $\ell ^{\infty }-$ interpolating
 sequence $K=(k_{j})_{j}$ such that $TR_{K}=I_{Z}$ (and $V_{K}$
  is bounded below) ; \par 
\quad (ii) the space ${\mathcal{E}}$ is weak*-complemented in $Z$ and
 $X_{K}$ is (norm-) complemented in $X.$\par 
\quad b) Conversely if $K=(k_{j})_{j}$ is an $\ell ^{\infty }-$ interpolating
 sequence such that $X_{K}$ is complemented in $X,$ say $X=X_{K}\oplus
 Y$  then $Y^{\perp }$ is weak*-homeomorphic to $\ell ^{\infty }.$
\end{Thrm}
\quad Of course the case $X=L^{1}/H_{0}^{1}$ gives the result announced
 by the title of this subsection.\ \par 
\quad Proof.\ \par 
\quad a) By Proposition~\ref{FB0} we can view $T$ as a weak*-continuous,
 bounded below map from $\ell ^{\infty }$ into $Z$ with range
 ${\mathcal{E}}.$ Thus $T=S^{*}$ and $S(\in {\mathcal{L}}(X,\ell
 ^{1}))$ is onto (by virtue of the classical duality result recalled
 in the introduction). Hence there exists a sequence $K=(k_{n})_{n}$
 in $X$ such that $S(k_{n})=\epsilon _{n}\ ,\ n\in {\mathbb{N}}.$
 Note that such a sequence is necesssarily biorthogonal to the
 sequence $(g_{n})_{n}=T(\epsilon _{n})_{n}$ in $Z$ ; indeed \ \par 
\quad \quad \quad $\delta _{j,n}={\left\langle{\epsilon _{j},\epsilon _{n}}\right\rangle}={\left\langle{S(k_{j}),\epsilon
 _{n}}\right\rangle}={\left\langle{k_{j},T\epsilon _{n}}\right\rangle}={\left\langle{k_{j},g_{n}}\right\rangle}.$\
 \par 
\quad Moreover, by the open mapping theorem, we may (and do) choose
 the sequence $K=(k_{n})_{n}$ to be bounded. Then $k_{n}=V_{K}(\epsilon
 _{n})$  for all $n,$ and hence $SV_{K}=I_{l^{1}}$ ; consequently,
 $V_{K}S\ (=V_{K}SV_{K}S)$  is a bounded projection in $X$ whose
 range $X_{K}$  is thus norm-complemented in $X.$ Similarly,
 $R_{K}T\ (=V_{K}^{*}S^{*}=(SV_{K})^{*}=I_{\ell ^{\infty }})$
 and $TV_{K}$ is a weak*-continuous projection with range ${\mathcal{E}}$
  which is therefore weak*-complemented in $Z.$\ \par 
\quad b) Here, since $K$ is an $\displaystyle \ell ^{\infty }$-IS,
 $V_{K}$  implements an isomorphism between $\ell ^{1}$ and $X_{K}$
 ; also we know, from the topological direct sum $\displaystyle
 X=X_{K}\oplus Y,$ that $X_{K}$ is isomorphic to the quotient
 of $X/Y.$ Therefore, by standard duality facts, the dual of
 $X_{K}$ is isomorphic to $Y^{\perp }.$\ \par 
Putting all this together : $\ell ^{\infty }=(\ell ^{1})^{*}\cong
 X_{K}^{*}\simeq Y^{\perp }$  we obtain the desired result. $\displaystyle
 \hfill\blacksquare $\ \par 
\ \par 
{\bf Remarks}\ \par 
\quad 1. We come back briefly to the context of theorem~\ref{F7} to
 point out that in this case the subspace ${\mathcal{E}}$ is
 also norm-complemented in the norm-closed subalgebra $\displaystyle
 {\mathcal{A}}_{\Omega }$ consisting of all the functions in
 $\displaystyle H^{\infty }$ continuously extendable to $\displaystyle
 {\mathbb{D}}\cup \Omega .$ Indeed the map ${\mathcal{J}}$ (associated
 to $\displaystyle (a_{n})_{n}$) makes sense on $\displaystyle
 {\mathcal{A}}_{\Omega }$ and so does the (norm-bounded) projection
 $\displaystyle P.$\ \par 
\quad 2. The map $T$ (from $\ell ^{\infty }$ into $H^{\infty }$) implementing
 the weak*-homeomorphism is usually called the (bounded) linear
 extension map (cf.~\cite{BiF10} , Chap. VII, Section 2).\ \par 
\quad Thus in the more general setting above we will say that the sequence
 $K=(k_{j})_{j}$ in $X$ is a boundedly linearly extendable sequence
 (notation: $BLES$) in $X$ for $Z$ if $R_{K}$ admits a right
 inverse (i.e. there exists $T\in {\mathcal{L}}(\ell ^{\infty
 },Z)$ such that $RT=I_{\ell ^{\infty }}$) (or equivalently if
 $V_{K}$ admits a left inverse).\ \par 
\ \par 
\quad It follows immediately from these definitions and characterizations
 that an $\displaystyle \ell ^{\infty }$-IS $K$ such that $X_{K}=X$
 is a $BLES.$ Thus (via standard facts about duality of subspaces
 and quotient spaces)  a sequence $K$ is an $\displaystyle \ell
 ^{\infty }$-IS in $X$ for $Z=X^{*}$ if and only if it is a $BLES$
 in $X_{K}$ for $\tilde Z=(X_{K})^{*}=Z/(^{\perp }X_{K}).$\ \par 
\quad One more important fact coming out right away from the above
 is that there exists an $\ell ^{\infty }-IS$ in $X$ for $Z=X^{*}$
  if and only if $X$ contains (isomorphically) a copy of $\ell
 ^{1},$ that is, in turn by virtue of the celebrated Rosenthal's
 $\ell ^{1}-$ theorem (cf.~\cite{BiF21} or [ Chap.1] in~\cite{BiF22}),
 if and only if there exists in $X$ a bounded sequence $(x_{n})_{n}$
 which admits no weakly Cauchy subsequence. Recall that a sequence
 $(y_{n})_{n}$ is weakly Cauchy if for any $h$ in $X^{*}$ the
 numerical sequence $({\left\langle{y_{n},h}\right\rangle})_{n}$
  converges. Rosenthal's $\ell ^{1}-$ theorem says more precisely
 that a bounded sequence without weakly Cauchy subsequence admits
 a subsequence equivalent to the canonical basis of $\ell ^{1}.$
 In particular we deduce from the above that there are such sequences
 in $L^{1}/H_{0}^{1}.$\ \par 
\quad 3. Theorem 10 in Chapter V of~\cite{BiF23} (which, as mentioned
 there, is essentially Theorem 4 of~\cite{BiF11} though the latter
 does not consider weak*-topology) states also (again, up to
 minor variations in the formulation) that Assertions a) and
 b) in Theorem~\ref{FB2} above are equivalent to the existence
 a copy of $c_{0}$ in $X^{*}.$ This is done via the (nice and
 rather deep) Bessaga-Pelcynski Selection principle and we refer
 to~\cite{BiF23}, Chap. V, for the details. We limit ourselves
 to pointing out that any bounded linear map $T\ :\ c_{0}\rightarrow
 X^{*}$  admits a bounded linear extension $\tilde T\ :\ \ell
 ^{\infty }\rightarrow X^{*}$  (just take $\tilde T=S^{*}$ with
 $S=(T^{*})_{\mid X}$ ). (Of course, there is no reason why $T$
 bounded below would imply $\tilde T$  bounded below.)\ \par 

\section{About $\ell ^{p}-$ interpolating sequences ($1\leq
 p<\infty $).~\label{F9}}
\ \par 

\subsection{The case $p=1$ ; Sidon sets.}
\ \par 
\quad From a theoretical point of view (and in terms of norm-isomorphism)
  the question of the existence of copies of $\ell ^{1}$ in a
 given Banach space $Y$ is settled by the celebrated Rosenthal's
 $\ell ^{1}-$ Theorem recalled above.\ \par 
\quad Concretely, we need a bounded sequence $G=(g_{j})_{j}$ in $Y$
 such that the associated map $V_{G}\ :\ \lambda =(\lambda _{n})_{n}\in
 \ell ^{1}\rightarrow \Sigma _{j}\lambda _{j}g_{j}$  is bounded below.\ \par 
\quad The Sidon sets provide a whole bunch of them "practically for
 free" in the disk algebra.\ \par 
\quad Indeed, {\sl one definition of a Sidon set of integers} is the following:\ \par 
{\sl       An infinite set }$\Gamma \subset {\mathbb{Z}}$ {\sl
 is a Sidon set (here, for }$n\in {\mathbb{Z}},\ \epsilon _{n}\
 ${\sl  is the function }$t\rightarrow e^{int}\ ${\sl  also identified
 as a function on the unit circle }${\mathbb{T}}${\sl ) if there
 exists a constant }$M$ {\sl such that, for any }$\alpha \in
 c_{00}(\Gamma ),\ \sum_{j}{\left\vert{\alpha _{j}}\right\vert
 }\leq M{\left\Vert{\sum_{j}{\alpha _{j}\epsilon _{j}}}\right\Vert}_{\infty
 }${\sl  (where }$\displaystyle c_{00}(\Gamma )${\sl  denotes
 the set of sequences with finite support in }$\Gamma ${\sl ).}\ \par 
\quad We have the following result.\ \par 

\begin{Thrm}
Let $\Gamma \subset {\mathbb{N}}$ a Sidon set ; then the (closed
 linear)  subspace ${\mathcal{E}}_{\Gamma }$ generated by the
 $\epsilon _{n},\ n\in \Gamma $  in the disc algebra $A({\mathbb{D}})$
 is isomorphic to $\ell ^{1}$ ; moreover this isomorphism is
 a weak*-homeomorphism (considering ${\mathcal{E}}_{\Gamma }$
 as a subspace of $H^{\infty }$).
\end{Thrm}
\quad Proof.\ \par 
The fact that $\Gamma $ is a Sidon set means (with the above
 definition)  that the (clearly valued in $A({\mathbb{D}})$ since
 here $\Gamma \subset {\mathbb{N}}$)  linear map $V_{\Gamma }\
 :\ \alpha \rightarrow \sum_{j\in \Gamma }{\alpha _{j}\epsilon
 _{j}}$ is well-defined, bounded and bounded below. This proves
 the first assertion. To conclude one checks that $V_{\Gamma
 },$ seen as a map from $\ell ^{1}$ in $H^{\infty },$ is the
 adjoint of the map $S$  defined, for $\lbrack f\rbrack \in L^{1}/H_{0}^{1},$
 by $S(\lbrack f\rbrack =(c_{-n}(f))_{n\in \Gamma }.$ Indeed
 $S$ is well-defined since if $\lbrack f\rbrack =\lbrack g\rbrack
 $ then $(f-g\mathrm{ }\mathrm{)}\mathrm{ }\in H_{0}^{1}$  and
 consequently $c_{-n}(f-g)\ (=\int{(f-g)\epsilon }_{n})=0$ for
 all $n\in {\mathbb{N}}$ and, in particular, for all $n\in \Gamma
 $ ; in addition $S$ is clearly linear, bounded and valued in
 $c_{0}(\Gamma ).$ Thus $V_{\Gamma }$ is weak*-continuous and
 implements a weak*-homeomorphism from $\ell ^{1}(\Gamma )$ onto
 ${\mathcal{E}}_{\Gamma }=V_{\Gamma }(\ell ^{1}(\Gamma )).$\ \par 
\ \par 

\begin{Exmp}
Among the Sidon sets contained in ${\mathbb{N}}$ we have the
 (set of values) of $q-$ lacunary sequences, i.e., sets of the
 form $\lbrace \lambda _{k};\ k\in {\mathbb{N}}\rbrace $  where
 the integers $\lambda _{k}$ satisfy $\ \frac{\lambda _{k+1}}{\lambda
 _{k}}\geq q$  for some (fixed) $q>1.$
\end{Exmp}
\ \par 

\begin{Rmrq}
In view of the result in Subsection~\ref{F13}, it is natural
 to ask if one can obtain a weak*-copy, say ${\mathcal{E}},$
 of $\ell ^{1}\ (=(c_{0})^{*})$  in $H^{\infty }\ (=(L^{1}/H_{0}^{1})^{*})$
 with, in addition, ${\mathcal{E}}$ weak*-complemented in $H^{\infty
 }.$ Following the same general and standard duality argument
 as in the proof of Theorem~\ref{FB2}, if this were the case,
 one would obtain easily that the subspace $\ ^{\perp }{\mathcal{E}}$
 has a complement (in $L^{1}/H_{0}^{1}$)  which is isomorphic
 to $c_{0}$ (everything so far is valid if we consider, instead
 of $H^{\infty }\ (=(L^{1}/H_{0}^{1})^{*}),$ an arbitrary dual
 Banach space $Z=X^{*}$) ; but the space $L^{1}/H_{0}^{1},$ being
 of cotype $2$ (cf.~\cite{BiF26}) , cannot contain a copy of
 $c_{0}.$ Thus we cannot have a weak*, weak*-complemented copy
 of $\ell ^{1}$ in $L^{1}/H_{0}^{1}.$
\end{Rmrq}
\ \par 

\subsection{The case $1<p<\infty .$}
\ \par 
\quad In this section $X$ is a separable reflexive Banach space, $Z=X^{*}$
  and $Z^{*}$ is identified with $X.$ We denote by $p'$ the conjugate
 of $p\ (\frac{1}{p}+\frac{1}{p'}=1).$\ \par 
\quad The following definitions are natural adaptations of the $p=\infty
 $  case.\ \par 

\begin{Dfnt}
A sequence $K=(k_{j})_{j}$ in $X$ is called: \par 
\quad a) an $\ell ^{p}-$ interpolating sequence for $Z$ ($\ell ^{p}-IS$
  for short) if the range of the linear map $R(=R_{K})\ :\ g\in
 Z\ \rightarrow \ ({\left\langle{k_{j},g}\right\rangle})_{j}\in
 \ell ^{0}$ contains $\ell ^{p}$ ;\par 
\quad b) an $\ell ^{p}-$ boundedly linearly extendable sequence for
 $Z$ ($\ell ^{p}-BLES$ for short) if there exists $T\in {\mathcal{L}}(\ell
 ^{p},Z)$  such that $\forall \lambda =(\lambda _{j})_{j}\in
 \ell ^{p},\ \ R\circ T(\lambda )=\lambda .$
\end{Dfnt}
\ \par 
\quad We list a few elementary observations and some examples before going on.\ \par 
\quad 1. An $\ell ^{p}-BLES$ is an $\ell ^{p}-IS.$\ \par 
\ \par 
\quad 2. If $K=(k_{j})_{j}$ is an $\ell ^{p}-BLES$ via $T\in {\mathcal{L}}(\ell
 ^{p},\ Z)$  then the sequence $(g_{j})_{j}$ defined by $g_{j}=T\epsilon
 _{j}$  for $j\in {\mathbb{N}},$ is biorthogonal to the sequence
 $K$ and is a so-called $p-$ Hilbertian sequence, that is, there
 exists a constant $C$ such that\ \par 
\quad \quad \quad $\ {\left\Vert{\sum_{j}{\alpha _{j}g_{j}}}\right\Vert}\leq C{\left\Vert{\alpha
 }\right\Vert}_{p},\ \alpha =(\alpha _{j})_{j}\in c_{00}\ \ \
 \ \ \ \ \ \ddag _{h}.$\ \par 
\quad Note also that, conversely, if $(g_{j})_{j}$ is a $p-$ Hilbertian
 sequence in $Z$ then there exists $T\in {\mathcal{L}}(\ell ^{p},Z)$
  such that $T\epsilon _{j}=g_{j}.$\ \par 
\ \par 
\quad 3. In case $X$ is an Hilbert space then of course $X_{K}$ is
 (orthogonally)  complemented, $Z$ may be identified with $X$
 and $\tilde Z$ with $X_{K}$ itself, and (since any bounded linear
 map from a subspace can be boundedly linearly extended to the
 whole space) the notions of $IS$ and $BLES$ coincide.\ \par 
\quad Like in the previous case a fundamental class of examples is
 provided by standard interpolation associated to sequences in
 the open unit disk for the Hardy spaces $H^{p}.$\ \par 
\quad We briefly review this context:\ \par 
\quad The duality $H^{p}=(H^{p'})^{*}$ can be expressed by the mutual action\ \par 
\quad \quad \quad $\displaystyle \Phi _{g}(f)=:<f,g>=\int_{{\mathbb{T}}}{f\bar
 g\ }(:=\frac{1}{2\pi }\int_{0}^{2\pi }{f(e^{it})\bar g(e^{it})dt),}\
 f\in H^{p},\ g\in H^{p'}.$\ \par 
\quad Note that in this approach $g\rightarrow \Phi _{g}$ is a conjugate
 linear isomorphism (isometric if $p=2$) from $H^{p'}$ onto $(H^{p})^{*}.$\
 \par 
\quad Starting with a sequence $(z_{n})_{n}$ in ${\mathbb{D}}$ , we
 define the sequence $K=(k_{j})_{j}$ in $H^{p'}$ by\ \par 
\quad \quad \quad $\displaystyle k_{j}(z)=(1-\left\vert{z_{j}}\right\vert ^{2})^{1/p}\frac{1}{1-\bar
 z_{j}z}.$\ \par 
\quad Thus our map $R_{K}\ (:\ H^{p}\rightarrow \ell ^{0})$ is here
 given by $R_{K}(g)=((1-\left\vert{z_{j}}\right\vert ^{2})^{1/p}g(z_{j})$
 that is exactly the map $T_{p}$ introduced in the development
 of interpolation theory in~\cite{BiF24}.\ \par 
\quad The following proposition is essentially a reformulation of observation
 2 above.\ \par 

\begin{Prps}
A sequence $K=(k_{j})_{j}$ in $X$ is an $\ell ^{p}-BLES$ for
 $Z$ iff it admits a biorthogonal $p-$ Hilbertian sequence.
\end{Prps}
\quad Once we have a bounded linear map from $\ell ^{p}$ into $Z$ we
 just need this map to be bounded below to have an isomorphic
 copy of $\ell ^{p}$  in $Z$ or in other words (we need) the
 existence of a constant $c>0$  such that\ \par 
\quad \quad \quad $\displaystyle \ {\left\Vert{\sum_{j}{\alpha _{j}g_{j}}}\right\Vert}\geq
 c{\left\Vert{\alpha }\right\Vert}_{p},\ \alpha =(\alpha _{j})_{j}\in
 c_{00}\ \ \ \ \ \ \ \ \ \ddag _{b}$\ \par 
which can be also expressed by saying that the sequence $(g_{n})_{n}$
  is a $p-$ Besselian sequence.\ \par 
(We advert the reader that there is no universal agreement on
 the terminology regarding the notions of Hilbertian and Besselian
 sequences.)\ \par 
\quad A sequence $(g_{n})_{n}$ which is both $p-$ Besselian and $p-$
 Hilbertian is called a $p-$ Riesz sequence. In other words,
 a basis (in any Banach space) is a $p-$ Riesz basis iff it is
 equivalent to the canonical basis of $\ell ^{p}.$\ \par 
\quad Thus the above translates immediately in the following characterization
 of copies of $\ell ^{p}$ in a Banach space.\ \par 

\begin{Thrm}
A (closed) subspace ${\mathcal{M}}$ of $Z$ is isomorphic to $\ell
 ^{p}$  if and only if it admits a $p-$ Riesz basis.
\end{Thrm}
\quad The previous discussion might lead one to think that the existence
 of a $p-$ Hilbertian sequence $(g_{n})_{n}$ in $Z$ which in
 addition is basic (that is, is a basis for the closed subspace
 ${\mathcal{M}}$  generated by this sequence) is enough to make
 ${\mathcal{M}}$ is isomorphic to $\ell ^{p}.$\ \par 
\quad We now describe elementary examples showing that this is not the case.\ \par 

\begin{Prps}
~\label{FB3}( $1<p<\infty $ ) In the space $\ell ^{p}$ there exists\par 
\quad (i) a basis $(u_{n})_{n}$ which is $p-$ Hilbertian but non $p-$
 Besselian (hence non $p-$ Riesz) ; \par 
\quad (ii) a basis $(v_{n})_{n}$ which is $p-$ Besselian but non $p-$
 Hilbertian; \par 
\quad (iii) a basis which is neither $p-$ Hilbertian nor $p-$ Besselian. 
\end{Prps}
\ \par 
\quad The construction is based on the following finite dimensional
 computational lemma.\ \par 

\begin{Lmm}
($1<p<\infty $) In the space $({\mathbb{C}}^{N+1},{\left\Vert{\
 }\right\Vert}_{p})$  there exists a basis $(e_{0},..,e_{N})$
 such that, denoting $(f_{0},..,f_{N})$  its dual basis, we have: \par 
\quad a) it is, as well as its dual basis in $({\mathbb{C}}^{N+1},{\left\Vert{\
 }\right\Vert}_{p'})$ , norm bounded by $2$ ; \par 
\quad b) for $\lambda \in {\mathbb{C}}^{N+1}{\left\Vert{\sum_{j=0}^{N}{\lambda
 _{j}e_{j}}}\right\Vert}_{p}\leq 3{\left\Vert{\lambda }\right\Vert}_{p}$
  ; \par 
\quad c) there exists $\lambda \in {\mathbb{C}}^{N+1}$ such that $\
 {\left\Vert{\sum_{j=0}^{N}{\lambda _{j}e_{j}}}\right\Vert}_{p}<\rho
 _{N}{\left\Vert{\lambda }\right\Vert}_{p}$  with $\mathrm{l}\mathrm{i}\mathrm{m}_{N\rightarrow
 \infty }\rho _{N}=0.$
\end{Lmm}
\quad Proof (of lemma).\ \par 
Let $(\epsilon _{j})_{0\leq j\leq N}$ denote the canonical basis
 in ${\mathbb{C}}^{N+1}.$\ \par 
Setting\ \par 
\quad \quad \quad $e_{j}=\epsilon _{j}+\frac{1}{N+1}v,\ j=1,..,N,\ e_{0}=\rho v,\
 f_{0}=\rho 'v,\ f_{j}=\epsilon _{0}+\epsilon _{j}$\ \par 
with $\rho =(N+1)^{-1/p}$ and $\rho '=(N+1)^{-1/p'}$ ensures that\ \par 
\quad \quad \quad $\displaystyle 1={\left\langle{e_{0},f_{0}}\right\rangle}={\left\Vert{e_{0}}\right\Vert}_{p}={\left\Vert{f_{0}}\right\Vert}_{p'}=1$\
 \par 
and that $(e_{j})_{0\leq j\leq N},\ (f_{j})_{0\leq j\leq N}$
 are dual bases in ${\mathbb{C}}^{N+1}$ (of course here, ${\mathbb{C}}^{N+1}$
  is -algebraically- identified with its dual via the canonical
 bilinear form $\ {\left\langle{a,b}\right\rangle}=\sum_{j=0}^{N}{a_{j}b_{j}}$)\!\!\!\!
 .\ \par 
\quad Now, the inequalities in a) are obvious for the $f_{j}$'s and
 $e_{0}$  and follow from the following ones for the $e_{j}$'s
 ($j=1,...,N$)  :\ \par 
\quad \quad \quad $\displaystyle \ {\left\Vert{e_{j}}\right\Vert}_{p}\leq 1+{\left\Vert{\frac{1}{N+1}v}\right\Vert}_{p}=1+(N+1)^{-1/p'}\leq
 2.$\ \par 
\quad For $\lambda \in {\mathbb{C}}^{N+1}$ we have\ \par 
\quad \quad \quad $\displaystyle \ \sum_{j=0}^{N}{\lambda _{j}e_{j}}=\sum_{j=1}^{N}{\lambda
 _{j}\epsilon _{j}}+\lambda _{0}e_{0}+\frac{1}{N+1}(\sum_{j=1}^{N}{\lambda
 _{j}})v.$\ \par 
\quad Since $\ {\left\Vert{\frac{1}{N+1}v}\right\Vert}_{p}=(N+1)^{-1/p'}$
 and (via Holder's inequality)  $\ \left\vert{\sum_{j=1}^{N}{\lambda
 _{j}}}\right\vert \leq N^{1/p'}{\left\Vert{\lambda }\right\Vert}_{p},$
 it follows easily (using the triangular inequality)\ \par 
\quad \quad \quad $\displaystyle \ {\left\Vert{\sum_{j=0}^{N}{\lambda _{j}e_{j}}}\right\Vert}_{p}\leq
 3{\left\Vert{\lambda }\right\Vert}_{p}.$\ \par 
Finally with regard to c) , with $\lambda =(0,1,1,..,1)$ we get\ \par 
\quad \quad \quad $\displaystyle \ {\left\Vert{\lambda }\right\Vert}_{p}=N^{1/p},\
 {\left\Vert{\sum_{j=0}^{N}{\lambda _{j}e_{j}}}\right\Vert}_{p}=\frac{1}{N+1}{\left\Vert{(N,1,...,\
 1)}\right\Vert}_{p}=\frac{1}{N+1}N^{1/p}(1+N^{p-1})^{1/p},$\ \par 
and hence $\displaystyle \rho _{N}:=\frac{{\left\Vert{\sum_{j=0}^{N}{\lambda
 _{j}e_{j}}}\right\Vert}_{p}}{{\left\Vert{\lambda }\right\Vert}_{p}}=\frac{(1+N^{p-1})^{1/p}}{1+N}\rightarrow
 0$ when $\displaystyle N\rightarrow \infty ,$ leading to the
 desired conclusion. $\displaystyle \hfill\blacksquare $\ \par 
\quad We now finish the proof of Proposition~\ref{FB3}.\ \par 
\quad Proof  (of proposition).\ \par 
Consider (with as usual $(\epsilon _{j})_{j\geq 0}$  the canonical
 basis in $\ell ^{p}$) the standard decomposition of $\ell ^{p}$
 into the ($\ell ^{p}$) direct sum of its $(N+1)$ dimensional
 $\ell _{N+1}^{p}$ subspaces:\ \par 
\quad \quad \quad $\ell _{2}^{p}=\bigvee_{k=0}^{1}{\epsilon _{k}},..\ ,\ell _{N+1}^{p}=\bigvee_{a_{N}\leq
 k<a_{N+1}}{\epsilon _{k}},$\ \par 
where $a_{1}=0$ and for $N>1$ $a_{N}\ (=a_{N-1}+N)=\frac{(N+2)(N-1)}{2}.$\ \par 
\quad To each of these subspaces we apply the above lemma and denote
 $(e_{k})_{a_{N}\leq k<a_{N+1}}$  the basis thus obtained and
 $(f_{k})_{a_{N}\leq k<a_{N+1}}$ the dual one.\ \par 
\quad It is clear that $(e_{k})_{k\geq 0}$ is a basis in $\ell ^{p}$
 which is $p-$ Hilbertian (indeed for all $\lambda \in c_{00},\
 {\left\Vert{\sum_{j}{\lambda _{j}e_{j}}}\right\Vert}_{p}\leq
 3{\left\Vert{\lambda }\right\Vert}_{p}$)  with a dual bounded
 basis (namely $(f_{k})_{k\geq 0}$ ) and this basis $(e_{k})_{k\geq
 0}$ cannot be $p-$ Besselian : indeed, from c) in the previous
 lemma ; in fact rather from its proof) we get for each $N$ an
 element $\lambda ^{(N)}$ (in $\ell _{N+1}^{p}$ hence in $c_{00}$)
  such that\ \par 
\quad \quad \quad $\ {\left\Vert{\lambda ^{(N)}}\right\Vert}_{p}=N^{1/p}\rightarrow
 \infty $  as $N\rightarrow \infty ,$\ \par 
while $\ {\left\Vert{\sum_{j}{\lambda _{j}^{(N)}}}\right\Vert}_{p}$
  is bounded.\ \par 
\quad Thus $(e_{k})_{k\geq 0}$ is an example for (i) and, by duality,
 $(f_{k})_{k\geq 0}$  is an example for (ii).\ \par 
As to (iii), we can just take $w_{k}=e_{k}$ for $a_{N}\leq k<a_{N+1}$
  when $N$ is even and $w_{k}=f_{k}$ for $a_{N}\leq k<a_{N+1}$
 when $N$ is odd.\ \par 
\ \par 

\begin{Rmrq}
The first example of a Hilbertian basis not Besselian was provided
 by Babenko~\cite{BiF9} in a Hilbert space context ; the above
 one (again in the Hilbert space setting) comes from the first
 named author's thesis~\cite{BiF25}  but was never published elsewhere.
\end{Rmrq}
\ \par 
\ \par 
\ \par

\end{document}